\documentclass[smallextended,envcountsect,envcountsame,envcountreset]{svjour3}
\usepackage{mathptmx}    
\usepackage{amsmath}
\usepackage{amssymb}
\usepackage{xspace}
\usepackage{gastex}

\smartqed

\newenvironment{Example}{\begin{example}}{\end{example}}
\newenvironment{Proof}{\begin{proof}}{\qed\end{proof}}
\newenvironment{Remark}{\begin{remark}}{\end{remark}}

\newcommand{\field}{\mathbb{F}}
\newcommand{\gf}[1]{GF({#1})}
\newcommand{\sdif}{\mathrm{\Delta}}
\newcommand{\invm}{\mathrm{inv}}
\newcommand{\vertexrem}{\setminus}
\newcommand{\dmatroid}{$\Delta$-matroid\xspace}
\newcommand{\dmatroids}{$\Delta$-matroids\xspace}
\newcommand{\dual}{\mathop{\bar{*}}}

\renewcommand{\emptyset}{\varnothing}

\begin{document}

\title{Quaternary Bicycle Matroids and the Penrose Polynomial for Delta-Matroids}
\titlerunning{Quaternary Bicycle Matroids and the Penrose Polynomial for Delta-Matroids}

\author{Robert Brijder
\thanks{R.B.\ is a postdoctoral fellow of the Research Foundation -- Flanders (FWO).}
 \and Hendrik Jan Hoogeboom}
\institute{R. Brijder \at
Hasselt University and Transnational University of Limburg, Belgium\\
\email{robert.brijder@uhasselt.be}
\and
H.J. Hoogeboom \at
Leiden Institute of Advanced Computer Science, Leiden University, The Netherlands}

\date{Received: date / Accepted: date}

\maketitle

\begin{abstract}
In contrast to matroids, vf-safe delta-matroids have three kinds of minors and are closed under the operations of twist and loop complementation. We show that the delta-matroids representable over $\gf{4}$ with respect to the nontrivial automorphism of $\gf{4}$ form a subclass of the vf-safe delta-matroids closed under twist and loop complementation. In particular, quaternary matroids are vf-safe.

Using this result, we show that the matroid of a bicycle space of a quaternary matroid $M$ is obtained from $M$ by using loop complementation. As a consequence, the matroid of a bicycle space of a quaternary matroid $M$ is independent of the chosen representation. This also leads to, e.g., an extension of a known parity-type characterization of the bicycle dimension, a generalization of the tripartition of Rosenstiehl and Read [Ann.\ Disc.\ Math.\ (1978)], and a suitable generalization of the dual notions of bipartite and Eulerian binary matroids to a vf-safe delta-matroids.

Finally, we generalize a number of results concerning the Penrose polynomial from binary matroids to vf-safe delta-matroids. In this general setting the Penrose polynomial turns out to have a recursive relation much like the recursive relation of the Tutte polynomial.
\end{abstract}

\section{Introduction}
It turns out that various results related to matroids are more generally and more efficiently obtained in the more general context of delta-matroids (or \dmatroids) defined by Bouchet~\cite{bouchet1987}. While matroids are closed under taking its dual, \dmatroids are closed under the more general operation of twist which can be viewed as a ``partial dual''. In \cite{BH/PivotLoopCompl/09} another operation for \dmatroids was defined, called loop complementation. Since \dmatroids in general are not closed under loop complementation, when studying loop complementation we often restrict to the (minor-closed) class of \dmatroids closed under both twist and loop complementation, called vf-safe \dmatroids. As twist and loop complementation form a group, vf-safe \dmatroids enjoy various interesting properties. In particular, vf-safe \dmatroids allow for three kinds of minors (a third in addition to the usual deletion and contraction). It was shown in \cite{BH/NullityLoopcGraphsDM/10} that binary \dmatroids (which includes all binary matroids) are vf-safe. In fact, the class of binary \dmatroids is closed under both twist and loop complementation. We generalize this result to $\gf{4}$. First we consider an extension of the notion of representability of \dmatroids over some field $\field$ using the notion of $\alpha$-symmetric matrices, where $\alpha$ is an automorphism of $\field$. While quaternary \dmatroids are not always vf-safe, we show that the \dmatroids representable over $\gf{4}$ \emph{with respect to} the nontrivial automorphism $\invm$ of $\gf{4}$ form a subclass of the vf-safe \dmatroids closed under twist and loop complementation, cf.\ Theorem~\ref{thm:inv_repr_gf4_vfclosed}. As a consequence, every quaternary matroid is vf-safe, cf.\ Corollary~\ref{cor:quat_vfsafe}.

Next we consider the effect of loop complementation on quaternary matroids (and binary matroids in particular). A ``full'' twist applied to a matroid obtains its dual matroid. We show in this paper that a ``full'' loop complementation applied to a quaternary matroid $M$ followed by taking all maximal sets results in a matroid representing \emph{any} bicycle space corresponding to $M$, cf.\ Theorem~\ref{thm:bicyclemat_char}. As a consequence, the matroids of two bicycle spaces corresponding to a quaternary matroid $M$ are equal. More generally, we show that loop complementation on a subset $Y$ of the ground set corresponds to the matroid of any bicycle space corresponding to $M$ \emph{relative to} $Y$ (this extended notion of bicycle space is defined in \cite{AignerPenroseBinMat} for binary matroids). The link between loop complementation and bicycle spaces has various consequences and explains why bicycle spaces appear often and in unexpected ways in the literature. We consider a number of results concerning the bicycle space of binary matroids and extend them to vf-safe \dmatroids. For example, we show that the well-known principal tripartition result of Rosenstiehl and Read \cite{Rosenstiehl1978195} can be generalized to vf-safe \dmatroids. Also, the notions of Eulerian matroid and bipartite matroid are dual for binary matroids. We show that these notions can be linked to loop complementation, and this link suggests alternative definitions of Eulerian and bipartite matroid that coincide for binary matroids but are (unlike the usual definitions) dual for the larger class of vf-safe matroids (or, indeed, vf-safe \dmatroids).

The final application of the results concerning the link between loop complementation and bicycle spaces is the Penrose polynomial. The Penrose polynomial is introduced by Penrose \cite{Penrose/Tensors/1971} to study the four-color conjecture (it was not yet a theorem then), see \cite{Aigner01PenroseGraphsMatroids} for a survey of the Penrose polynomial. The Penrose polynomial is defined in \cite{AignerPenroseBinMat} for binary matroids $M$ in general and in terms of the dimensions of the bicycle spaces of $M$ relative to the subsets $Y$ of the ground set of $M$. Using the obtained results concerning bicycle matroids, we straightforwardly generalize the Penrose polynomial to vf-safe \dmatroids, and then show that this polynomial allows for a recursive relation that \emph{characterizes} the polynomial. For the case of binary \dmatroids, we also formulate this polynomial as a graph polynomial with a recursive relation. As the class of binary matroids is not closed under loop complementation (in contract to the class of binary \dmatroids), we remark that this recursive relation is \emph{not} valid when restricting to the narrow viewpoint of binary matroids. This provides a further example of why the more general viewpoint of (vf-safe) \dmatroids is often worthwhile to consider. In fact, another example is provided by Chun et al.\ \cite{ChunMoffatt/DeltaM/EmbeddedGraphs}, where it is shown, using a preprint of this paper on arXiv, that the Penrose polynomial for vf-safe \dmatroids and its recursive relation turn out to generalize the Penrose polynomial for graphs embedded in surfaces of \cite{EllisMonaghan/PenroseEmb/2011} and its recursive relation. Finally, we consider evaluations of the Penrose polynomial for vf-safe \dmatroids inspired by the evaluations of the Penrose polynomial for binary matroids of \cite{AignerPenroseBinMat}.

\section{Preliminaries} \label{sec:prelim}
We first recall some basic notions and results.

\subsubsection*{Principal pivot transform.}
Let $V$ and $W$ be finite sets. We consider $V
\times W$-matrices $A$, i.e., matrices where the rows are
indexed by $V$ and the columns by $W$. The rows and
columns of $A$ are not ordered (note that matrix inversion,
rank, etc.\ are defined for such matrices). For $X \subseteq V$
and $Y \subseteq W$, the $X \times Y$-submatrix of $A$ is
denoted by $A[X,Y]$. We write simply $A[X]$ to denote $A[X,X]$. We define the \emph{deletion} of $X$ in $A$ by $A\vertexrem X = A[V\setminus X, W\setminus X]$.

Let $A$ be a $V \times V$-matrix (over an arbitrary field
$\field$), and let $X \subseteq V$ be such that the principal
submatrix $A[X]$ is nonsingular. The \emph{principal pivot
transform} (PPT) of $A$ on $X$, denoted by
$A*X$, is defined as follows \cite{tucker1960}. Let $ A =
\bordermatrix{ & \scriptstyle X & \scriptstyle V\setminus X \cr
\scriptstyle X & P & Q \cr \scriptstyle V\setminus X & R & S }
$, then
$
A*X = \bordermatrix{
& \scriptstyle X & \scriptstyle V\setminus X \cr
\scriptstyle X & P^{-1} & -P^{-1} Q \cr
\scriptstyle V\setminus X & R P^{-1} & S - R P^{-1} Q
}
$.
PPT has many applications and is well motivated as it can be viewed as a partial matrix inversion (full matrix inversion corresponds to the case $X=V$)~\cite{Tsatsomeros2000151}.

We now recall the following property of PPT.
\begin{proposition}[\cite{tucker1960,ParsonsTDP70}]\label{prop:tucker}
Let $A$ be a $V \times V$-matrix, and let $X\subseteq V$ be
such that $A[X]$ is nonsingular. Then, for all $Y \subseteq
V$, $\det((A*X)[Y]) = \det(A[X \sdif Y]) / \allowbreak \det(A[X])$. In particular, $(A*X)[Y]$ is nonsingular iff $A[X \sdif Y]$ is nonsingular.
\end{proposition}

It is easy to verify (by the above definition of PPT) that $-(A*X)^T = (-A^T)*X$ for all $X \subseteq V$ with $A[X]$ nonsingular.
As a consequence, if $A$ is skew-symmetric, i.e., $-A^T = A$ (we allow
nonzero diagonal entries in case $\field$ is of characteristic
two), then $A*X$ is skew-symmetric as well.

Hence, if $A$ is a $V \times V$-symmetric matrix over $\gf{2}$,
then so is $A*X$. We identify $V \times V$-symmetric matrices
$A$ over $\gf{2}$ with (undirected) graphs $G = (V,E)$ where $\{x,y\} \in E$
iff $A[\{x\},\{y\}] = 1$ (we allow $x=y$, i.e., loops); $A$ is called the \emph{adjacency matrix} of $G$. Hence,
we write, e.g., $G[X]$, $G\vertexrem X$, and $G*X$ to denote $A[X]$, $A\vertexrem X$, and $A*X$
respectively, where $A$ is the adjacency matrix of $G$. We will
often use simply $V$ to denote the vertex set of graph $G$
under consideration.

\subsubsection*{Twist and loop complementation on set systems.}
A \emph{set system} (over $V$) is a tuple $M = (V,D)$ with $V$
a finite set called the \emph{ground set} and $D \subseteq 2^V$
a family of subsets of $V$. Similar as with graphs, we will
often use $V$ to denote the ground set of set system $M$ under
consideration. We write simply $X \in M$ to denote $X \in D$. Set system $M$ is called \emph{proper} if $D \neq \emptyset$. We say that $M$ is \emph{equicardinal} if for all $X,Y \in M$, $|X| = |Y|$. A
set system is called \emph{even} if for all $X,Y \in M$, $|X|$ and $|Y|$ have equal parity.
We define, for $X \subseteq V$, the \emph{restriction}
of $X$ in $M$ by $M[X] = (X,D')$ where $D' = \{Y \in D
\mid Y \subseteq X\}$, and we define $M \vertexrem X = M [V \setminus X]$.
We define, for $X \subseteq V$, the \emph{twist} \cite{bouchet1987} of $M$ on $X$,
denoted by $M * X$, as $(V,D
* X)$, where $D * X = \{Y \sdif X \mid Y \in D\}$ and $\sdif$ denotes symmetric difference.
We say that $u \in V$ is a \emph{loop} in $M$ if none of the $X \in M$ contains $u$. We say that $u$ is a \emph{coloop} in $M$ if $M*V$ is a loop in $M$, i.e., all $X \in M$ contain $u$.
Moreover, we define, for $X \subseteq V$, \emph{loop
complementation} of $M$ on $X$, denoted by $M + X$, as
$(V,D')$, where $Y \in D'$ iff $|\{ Z \in M \mid Y \setminus X
\subseteq Z \subseteq Y \}|$ is odd
\cite{BH/PivotLoopCompl/09}. As recalled below, loop complementation for set systems turns out to generalize loop complementation for graphs.

We assume left associativity of set system
operations. Therefore, e.g., $M+X \vertexrem Y$ denotes
$(M+X) \vertexrem Y$.  Twist and loop complementation are
involutions and they commute on distinct elements. Hence for
$X,Y \subseteq V$ we have, e.g., $M * X * Y = M * (X \sdif Y)$,
$M + X + Y = M+ (X \sdif Y)$, and if $X\cap Y = \emptyset$,
$M * X + Y = M + Y * X$. It turns out that ${}*X$ and ${}+X$
(for any $X \subseteq V$) generate the group $S_3$ of
permutations on $3$ elements \cite{BH/PivotLoopCompl/09}. We
denote by ${}\dual X$ the third element ${}*X+X*X = {}+X *X +X$ of order $2$,
called the \emph{dual pivot} on $X$. We have that
$M \dual X$ is equal to $(V,D')$, where $Y \in D'$ iff $|\{ Z
\in M \mid Y \subseteq Z \subseteq Y \cup X \}|$ is odd. Let
$\min(M) = (V,\min(D))$ and $\max(M) = (V,\max(D))$, where
$\min(D)$ ($\max(D)$, resp.) are the sets in $D$ which are
minimal (maximal, resp.) with respect to set inclusion. We denote by
$d_M = \min_{Y\in M}(|Y|)$ the smallest cardinality among the
sets in $M$. It is observed in \cite{BH/PivotLoopCompl/09} that
$\min(M) = \min(M+X)$, thus $d_M = d_{M+X}$. Since $\min(M) =
\max(M*V)*V$, we have
similarly $\max(M) = \max(M\dual X)$.  
In case of singletons $\{u\}$, we also write $M \vertexrem
u$, $M*u$, etc.\ to denote $M \vertexrem \{u\}$, $M * \{u\}$,
etc.

\subsubsection*{Matroids and \dmatroids.}
We assume the reader is familiar with the basic notions
concerning matroids, which can be found, e.g., in
\cite{Welsh/MatroidBook,Oxley/MatroidBook-2nd}. We recall here \dmatroids and their link with matroids.

A proper set system $M$ is a \emph{\dmatroid
}~\cite{bouchet1987} if for all $X,Y \in M$ and for all $x \in X
\sdif Y$, there is a $y \in X \sdif Y$ ($x = y$ is allowed),
such that $X \sdif \{x,y\} \in M$. If $M$ is a \dmatroid, then
so is $M*X$ for all $X \subseteq V$.

A set system $M$ is an equicardinal \dmatroid iff $M$ is a matroid
described by its bases \cite{DBLP:conf/ipco/Bouchet95}.
Also, if $M$ is a \dmatroid, then $\min(M)$ and $\max(M)$ are matroids \cite{Bouchet_1991_67}.
In this paper, we assume, unless stated otherwise, that \emph{matroids are described by their bases}.
Hence, if $M$ is a matroid, then $M*V = M^*$ is the dual matroid of
$M$. The nullity and rank of a matroid $M$ are denoted by $\nu(M)$
and $\rho(M)$, respectively. Note that $\rho(M) = d_M$ and
$\nu(M) = d_{M*V}$ for matroids $M$, and $\rho(\min(M)) = d_{M}$ and $\nu(\max(M)) = d_{M*V}$ for \dmatroids $M$. Also note that the notions of loop and coloop for set systems coincide with the eponymous notions for matroids.

A \emph{minor} of a \dmatroid is a \dmatroid obtained by applying a (possibly empty) sequence of operations of the form ${}\vertexrem u$ and ${}*u\vertexrem u$ with $u \in V$. Note that the usual (matroid-theoretical) notion of deletion of $u$ for a matroid $M$ is equal to $M \vertexrem u$ if $u$ is not a coloop and equal to $M*u\vertexrem u$ otherwise. Similarly, contraction of $u$ for a matroid $M$ is equal to $M*u \vertexrem u$ if $u$ is not a loop and equal to $M \vertexrem u$ otherwise. Since, in this paper, we never apply $\vertexrem u$ to a coloop $u$ or ${}*u \vertexrem u$ to a loop $u$, the reader may think of ${}\vertexrem u$ and ${}*u\vertexrem u$ as the usual matroid-theoretical notions of deletion and contraction, respectively. A \dmatroid $M$ is called \emph{vf-safe} if any set system in the orbit of $M$ under ${}+{}$ and ${}*{}$ is a \dmatroid. The class of vf-safe \dmatroids is minor closed \cite{BH/NullityLoopcGraphsDM/10}. There are (delta-)matroids that are not vf-safe, such as the $6$-point line $U_{2,6}$, $P_6$, and the non-Fano matroid $F_7^-$, see \cite{Oxley/MatroidBook-2nd} for a description of these matroids. In fact, they are excluded minors for the class of vf-safe \dmatroids \cite{BH/NullityLoopcGraphsDM/10}.

\section{$\alpha$-symmetry and delta-matroids} \label{sec:alpha_symm_dmatroids}
Let $\alpha$ be an automorphism of a field $\field$. By abuse of notation, we extend $\alpha$ point-wise to vectors, matrices, and subspaces over $\field$. Hence for a $V \times V$-matrix $A = (a_{i,j})_{i,j \in V}$, we let $\alpha(A) = (\alpha(a_{i,j}))_{i,j \in V}$. Moreover, for subspace $L \subseteq \field^V$, we let $\alpha(L) = \{ \alpha(v) \mid v \in L\}$.

Let $A$ be a $V \times V$-matrix over some field $\field$, and let $\alpha$ be an involutive automorphism of $\field$, i.e., $\alpha^2(x) = x$ for all $x$ of $\field$ (the identity automorphism is considered an involutive automorphism here). Then $A$ is called \emph{$\alpha$-symmetric} if $\alpha(-A^T) = A$.
Note that if $A$ is $\alpha$-symmetric, then $A^T$ is $\alpha$-symmetric. Also note that $A$ is $\mathrm{id}$-symmetric with $\mathrm{id}$ the identity automorphism iff $A$ is skew-symmetric.

\begin{lemma} \label{lem:automorph_PPT}
Let $A$ be a $V \times V$-matrix over some field $\field$, and let $\alpha$ be an automorphism of $\field$. If $X \subseteq V$ is such that $A[X]$ is nonsingular, then $\alpha(A*X) = \alpha(A)*X$.
\end{lemma}
\begin{Proof}
For any nonsingular matrix $P$, we have $\alpha(P^{-1}) = \alpha(P)^{-1}$. Thus,
if $A
={}$ $ \bordermatrix{
& \scriptstyle X & \scriptstyle V\setminus X \cr
\scriptstyle X & P & Q \cr
\scriptstyle V\setminus X & R & S
}$, then in both cases we obtain:
$$
\bordermatrix{
& \scriptstyle X & \scriptstyle V\setminus X \cr
\scriptstyle X & \alpha(P)^{-1} & -\alpha(P)^{-1} \alpha(Q) \cr
\scriptstyle V\setminus X & \alpha(R) \alpha(P)^{-1} & \alpha(S) - \alpha(R) \alpha(P)^{-1} \alpha(Q)
}.
$$
\end{Proof}

If $A$ is $\alpha$-symmetric and $X \subseteq V$ is such that $A[X]$ is nonsingular, then $A*X$ is $\alpha$-symmetric. Indeed, $\alpha(-(A*X)^T) = \alpha((-A^T)*X) = \alpha(-A^T)*X = A*X$, where in the second equality we use Lemma~\ref{lem:automorph_PPT}.

Let $A$ be a $V \times V$ matrix. We define the set system
$\mathcal{M}_A = (V,D)$ where $X \in D$ iff $A[X]$ is
nonsingular. By convention, $A[\emptyset]$ is nonsingular. The next result is a straightforward extension of a result of \cite{bouchet1987} (the original formulation restricts to the case $\alpha = \mathrm{id}$).
\begin{lemma} [\cite{bouchet1987}]  \label{lem:alpha_symm_dmatroid}
Let $\alpha$ be an involutive automorphism of some field $\field$, and let $A$ be a $\alpha$-symmetric $V \times V$-matrix over $\field$. Then $\mathcal{M}_{A}$ is a \dmatroid.
\end{lemma}
\begin{Proof}
Let $X,Y \in \mathcal{M}_{A}$ and $x \in X \sdif Y$. If entry $A*X[\{x\}]$ is nonzero, then by Proposition~\ref{prop:tucker}, $X \sdif \{x\} \in \mathcal{M}_{A}$ and we are done. Thus assume that $A*X[\{x\}]$ is zero. Since $A[Y]$ is nonsingular, $A*X[X\sdif Y]$ is nonsingular by Proposition~\ref{prop:tucker}. Hence there is a $y \in X \sdif Y$ with entry $A*X[\{x\},\{y\}]$ nonzero (note that $x \neq y$). Since $A*X$ is $\alpha$-symmetric, $A*X[\{x,y\}]$ is of the form
$$
\bordermatrix{&x&y\cr
x & 0 & t_1 \cr
y & \alpha(-t_1) & t_2
}
$$
for some $t_1 \in \field \setminus \{0\}$ and $t_2 \in \field$. Thus $A*X[\{x,y\}]$ is nonsingular and $X \sdif \{x,y\} \in \mathcal{M}_{A}$.
\end{Proof}

We say that a \dmatroid $M$ is \emph{$\alpha$-representable} over $\field$, if $M = \mathcal{M}_A*X$ for some $\alpha$-symmetric $V \times V$-matrix $A$ and $X \subseteq V$. In this way, the notion of \emph{representable} from \cite{bouchet1987} coincides with $\mathrm{id}$-representable.

For a graph $G$, $\mathcal{M}_G$ is even iff $G$ has no loops. It is easy to verify that for graphs $G$ and $G'$, $\mathcal{M}_G = \mathcal{M}_{G'}$ iff $G = G'$. In fact, $G$ is uniquely determined by $V$ and the sets of cardinality $1$ and $2$ of $\mathcal{M}_G$, see \cite[Property~3.1]{Bouchet_1991_67}.

A $V \times V$-matrix $A$ over $\field$ is called \emph{principally  unimodular} (PU, for short) if for all $Y \subseteq V$, $\det(A[Y]) \in \{0,1,-1\}$. Note that any $V \times V$-matrix over $\gf{2}$ or $GF(3)$ is principally unimodular.

We now consider the field $\gf{4}$. Let us denote the unique nontrivial automorphism of $\gf{4}$ by $\invm$. Note that $\invm(x) = x^{-1}$ for all $x \in \gf{4} \setminus \{0\}$, and thus $\invm$ is an involutive automorphism.

\begin{theorem} \label{thm:inv_symm_PU}
Let $A$ be a $\invm$-symmetric $V \times V$-matrix over $\gf{4}$. Then $A$ is a principally unimodular.
\end{theorem}
\begin{Proof}
Recall that $1 = -1$ in $\gf{4}$. We have $\det(A) = \det(\invm(-A^T)) = \invm(\det(-A^T)) \allowbreak = \invm(\det(A))$. Thus $\det(A) \in \{0,1\}$.
\end{Proof}

\begin{Remark}
The proof of Theorem~\ref{thm:inv_symm_PU} essentially uses that the field $\field$ under consideration is of characteristic $2$, i.e., $\field = GF(2^k)$ for some $k \geq 1$, \emph{and} that $\field$ has an involutive automorphism $\alpha$ with only trivial fixed points (the set of fixed points form $\gf{2}$). The automorphisms $\alpha$ of $GF(2^k)$ are of the form $x \mapsto x^{p^\ell}$, with $1 \le \ell \le k$, and $\alpha$ is an involution when either $k=1$ (and thus $\ell=1$) or both $k$ is even and $\ell = k/2$. Moreover, for $\ell = k/2$ and $k$ even, the corresponding automorphism $\alpha$ has only trivial fixed points iff $\ell=1$. Consequently, the proof of Theorem~\ref{thm:inv_symm_PU} only works for $\alpha = \invm$ and $\field = \gf{4}$ (and, of course, $\alpha = \mathrm{id}$ and $\field = \gf{2}$).
\end{Remark}

For $X \subseteq V$, we define $A+X$ to be the $V \times V$ matrix with $A+X[\{x\},\{y\}] = A[\{x\},\{y\}]+1$ if $x=y \in X$ and $A+X[\{x\},\{y\}] = A[\{x\},\{y\}]$ otherwise. The following result is a straightforward generalization of a result of \cite{BH/PivotLoopCompl/09} formulated for the case $\field = \gf{2}$.
\begin{proposition} [Theorem~8 of \cite{BH/PivotLoopCompl/09}] \label{prop:pu_loopc}
Let $A$ be a principally unimodular $V \times V$-matrix over a field $\field$ of characteristic $2$. Then, for all $X \subseteq V$, $\mathcal{M}_{A+X} = \mathcal{M}_A+X$.
\end{proposition}
\begin{Proof}
It suffices to show the result for $X = \{u\}$ with $u \in V$. By the definition of loop complementation, we need to show that, for $Z\subseteq V$, $Z \in \mathcal{M}_{A+u}$ iff (1) $Z \in \mathcal{M}_A$ when $u \notin Z$ and (2) exactly one of $Z$, $Z \setminus \{u\}$ is in $\mathcal{M}_A$ when $u \in Z$.

Let $Z\subseteq V$. First assume that $u\notin Z$. Then $A[Z] = (A+u)[Z]$, thus $\det A[Z] = \det (A+u) [Z]$ and so $Z \in \mathcal{M}_{A+u}$ iff $Z \in \mathcal{M}_A$. Now
assume that $u\in Z$, which implies that $A[Z]$ and $(A+u)[Z]$ differ
in exactly one position: $(u,u)$. We may compute determinants by
Laplace expansion over the $u$-column, and summing minors. As
$A[Z]$ and $(A+u)[Z]$ differ at only the matrix-element $(u,u)$,
these expansions differ only by the minor $\det
A[Z\setminus \{u\}]$. Thus $\det (A+u)[Z] = \det A[Z] + \det
A[Z\setminus \{u\}]$, and this computation is in $\gf{2}$ as $A$ is PU and $\field$ of characteristic $2$. Hence $Z \in \mathcal{M}_{A+u}$ iff exactly one of $Z$, $Z \setminus \{u\}$ is in $\mathcal{M}_A$.
\end{Proof}

Note that, for a graph $G$ (i.e., a symmetric matrix over $\gf{2}$), $G+X$ is obtained from $G$ by complementing the existence of loops for the vertices in $X$, hence the name loop complementation for the set systems operation $+X$.

A \dmatroid $M$ is said to be \emph{representable} over $\field$, if $M = \mathcal{M}_A*X$ for some skew-symmetric $V \times V$-matrix $A$ and some $X \subseteq V$. A \dmatroid is said to be \emph{binary} if it is representable over $\gf{2}$. The following result is an adaption of the proof of \cite[Theorem~8.2]{BH/NullityLoopcGraphsDM/10} where it is shown that the class of binary \dmatroids is closed under twist and loop complementation.
\begin{theorem} \label{thm:inv_repr_gf4_vfclosed}
The class of \dmatroids $\invm$-representable over $\gf{4}$ is closed under twist and loop complementation.
\end{theorem}
\begin{Proof}
Let $M$ be a \dmatroid $\invm$-representable over $\gf{4}$. Then $M = \mathcal{M}_A*X$ for some $\invm$-symmetric $V \times V$-matrix $A$ over $\gf{4}$ and $X \subseteq V$. Let
$\varphi$ be a sequence of twist and loop complementations over $V$.
Let $W \in \mathcal{M}_A*X \varphi$, and consider now $\varphi'
= *X\varphi*W$. By the $S_3^V$ group structure of $*$ and $+$,
$\varphi'$ can be put in the following normal form:
$\mathcal{M}_A \varphi' = \mathcal{M}_A+Z_1*Z_2+Z_3$ for some
$Z_1,Z_2,Z_3 \subseteq V$ with $Z_1 \subseteq Z_2$. By Theorem~\ref{thm:inv_symm_PU}, $A$ is PU. By Proposition~\ref{prop:pu_loopc},
$\mathcal{M}_A + Z_1 = \mathcal{M}_{A+Z_1}$. Thus
$\mathcal{M}_A+Z_1*Z_2+Z_3 = \mathcal{M}_{A+Z_1}*Z_2+Z_3$. By
construction $\emptyset \in \mathcal{M}_A \varphi'$. Hence we
have $\emptyset \in \mathcal{M}_{A+Z_1}*Z_2$. Therefore $Z_2
\in \mathcal{M}_{A+Z_1}$ and so $A+Z_1[Z_2]$ is nonsingular. Thus, $A+Z_1*Z_2$ is defined.
Consequently, $A' = A+Z_1*Z_2+Z_3$ is defined and
$\mathcal{M}_A \varphi' = \mathcal{M}_{A'}$. Hence $M\varphi =
\mathcal{M}_A*X\varphi = \mathcal{M}_{A'}*W$ and thus $\invm$-symmetric matrix $A'$
represents $M\varphi$. Therefore, $M\varphi$ a \dmatroid $\invm$-representable over $\gf{4}$.
\end{Proof}
Consequently, the class of vf-safe \dmatroids contains the class of $\invm$-rep\-re\-sentable \dmatroids over $\gf{4}$ (and therefore also the class of binary \dmatroids). In contrast with Theorem~\ref{thm:inv_repr_gf4_vfclosed}, it is shown in \cite{BH/NullityLoopcGraphsDM/10} that there are \dmatroids $\mathrm{id}$-representable over $\gf{4}$ that are \emph{not} vf-safe.

\section{Quaternary matroids}

Let $M = (V, \mathcal{B})$ be a matroid representable over $\field$, and described by its bases. Let $B$ be a standard representation of $M$ over $\field$. Then $B$ is equal to
$$
\bordermatrix{&\scriptstyle X&\scriptstyle V\setminus X\cr
\scriptstyle X & I & S
}
$$
for some $X \in \mathcal{B}$, where $I$ is the identity matrix of suitable size. Let $\alpha$ be an involutive automorphism of $\field$. We define $R(B,\alpha)$ to be the $\alpha$-symmetric $V \times V$-matrix
$$
\bordermatrix{&\scriptstyle X&\scriptstyle V\setminus X\cr
\scriptstyle X & 0 & S\cr
\scriptstyle V\setminus X & \alpha(-S^T) & 0
}.
$$

We now recall the following result of de~Frayseix~\cite{DBLP:journals/dm/Fraysseix81} and Bouchet~\cite{bouchet1987} which states that a matroid is representable in the classical matroid sense iff it is representable in the \dmatroid sense. Therefore, the class of matroids representable over some field $\field$ is a subclass of the the class of \dmatroids representable over $\field$.
\begin{proposition} [Theorem~4.4 of \cite{bouchet1987}] \label{prop:bouchet_twist_matroid}
Let $M$ be a matroid representable over $\field$, let $\alpha$ be an involutive automorphism of $\field$, and let $B$ be a $X \times V$-matrix over $\field$ that is a standard representation of $M$. Then $\mathcal{M}_A = M*X$ with $A = R(B,\alpha)$.
\end{proposition}

The formulation of Proposition~\ref{prop:bouchet_twist_matroid} is slightly more general than the original formulation in \cite{bouchet1987}, which assumes $\alpha = \mathrm{id}$. However, note that if $A = R(B,\alpha)$ and $A'= R(B,\mathrm{id})$, then $A[Y]$ is nonsingular iff $A'[Y]$ is nonsingular for all $Y \subseteq V$. Thus $\mathcal{M}_A = \mathcal{M}_{A'}$.

In case $\field = \gf{2}$, we can view $A$ from Proposition~\ref{prop:bouchet_twist_matroid} as (an adjacency matrix representation of) a $(X,V\setminus X)$-bipartite graph $G$. Graph $G$ is often called the \emph{fundamental graph} of $M$ with respect to the basis $X \in M$, consisting of all edges $\{u,v\}$ such that $X \sdif \{u,v\}$ is a basis. If $Y \in M$, then $M*Y = \mathcal{M}_G *X *Y = \mathcal{M}_G *(X \sdif Y) = \mathcal{M}_{G *(X \sdif Y)}$ where in the last equality we use that $X \sdif Y \in M*X = \mathcal{M}_G$ since $Y \in M$. Therefore, every fundamental graph of $M$ can be obtained from $G$ by applying PPT \cite[Section~2]{DBLP:journals/ejc/Bouchet01}.

Hence, by Proposition~\ref{prop:bouchet_twist_matroid}, a matroid $M$ is representable over $\field$ in the usual (matroid) sense iff $M$ is $\alpha$-representable for some involutive automorphism $\alpha$ of $\field$ (recall that we allow $\alpha = \mathrm{id}$) iff $M$ is $\alpha$-representable for all involutive automorphisms $\alpha$ of $\field$.
Therefore, choosing $\alpha = \mathrm{id}$ may not necessarily be the most convenient extension of the matroid notion of representability to \dmatroids. Indeed, in view of Theorem~\ref{thm:inv_repr_gf4_vfclosed} and the remark below it, we argue that over $\gf{4}$, $\invm$-representability is a more natural extension of the matroid notion of representability to \dmatroids.

In particular, every quaternary matroid is a \dmatroid $\invm$-representable over $\gf{4}$. Hence by Theorem~\ref{thm:inv_repr_gf4_vfclosed} we have the following result, which was conjectured in \cite{BH/NullityLoopcGraphsDM/10}.
\begin{corollary} \label{cor:quat_vfsafe}
Every quaternary matroid is vf-safe.
\end{corollary}

\section{Bicycle matroids}
Let $v \in \field^V$ be a vector. The \emph{support} of $v$ is the set $X \subseteq V$ such that the entries of $X$ in $v$ are nonzero and entries of $V\setminus X$ in $v$ are zero. Let $L$ be a subspace of $\field^V$. We denote by $M(L)$ the matroid with ground set $V$ such that for all $C \subseteq V$, $C$ is a circuit of $M(L)$ iff there is a $v \in L$ with support $C$ and $C$ is minimal with this property among the nonempty subsets of $V$. Note that for a $X \times V$-matrix $A$, the matroid $M(\ker(A))$ equals the column matroid of $A$, denoted by $M(A)$. The \emph{orthogonal complement} of $L$, denoted by $L^{\perp}$, is $\{v \in \field^V \mid \langle u,v \rangle = 0 \mbox{ for all } u \in L \}$ where $\langle u,v \rangle = \sum_{x \in V} u(x)v(x)$ for all $u,v \in \field^V$.

Consider now the case $\field = \gf{4}$. Inspired by terminology from \cite{DBLP:journals/jct/Vertigan98}, we call $L \cap \invm(L^{\perp})$ the \emph{bicycle space} of $L$, and denote it by $\mathrm{BC}_L$. We have $\mathrm{BC}_{L^\perp} = \invm(\mathrm{BC}_{L})$. The dimension of $\mathrm{BC}_L$ is called the \emph{bicycle dimension} of $L$. More generally, for $Y \subseteq V$ and vector $v$ over $V$, denote by $\pi_Y(v)$ the vector obtained from $v$ by setting all entries of $V\setminus Y$ to $0$. Then we call $\{v \in L \mid \pi_Y(v) \in \invm(L^{\perp})\}$ the \emph{bicycle space} of $L$ \emph{relative} to $Y$ and we denote it by $\mathrm{BC}_L(Y)$. Note that $\mathrm{BC}_L(\emptyset) = L$ and $\mathrm{BC}_L(V) = \mathrm{BC}_L$.

We now recall the notion of bicycle space of a binary matroid $M$. Let $M$ be a binary matroid over $V$ and let $\mathcal{CS}_M$ be the cycle space of $M$, i.e., the subspace of $\gf{2}^V$ generated by the circuits of $M$. The \emph{bicycle space} of $M$ \emph{relative to} $Y \subseteq V$ is defined as $\mathrm{BC}_M(Y) = \{C \in \mathcal{CS}_M \mid C \cap Y \in \mathcal{CS}_M^{\perp}\}$ (where vectors over $\gf{2}$ are identified by by their support), see \cite{AignerPenroseBinMat}. We (may) consider $\gf{2}^V$ as a subspace of $\gf{4}^V$. Observe that if $L \subseteq \gf{2}^V$, then $\mathrm{BC}_L = L \cap L^{\perp}$. It is well known that $\mathcal{CS}_M$ is equal to the null space $\ker(A)$ of any binary representation $A$ of $M$, see \cite[Proposition~9.2.2]{Oxley/MatroidBook-2nd} (in particular, $\dim(\mathcal{CS}_M) = \nu(M)$). Thus $\mathrm{BC}_{\ker(A)}(Y) = \mathrm{BC}_M(Y)$ for all $Y \subseteq V$ and so the definition of $\mathrm{BC}_L(Y)$ is consistent with the definition of $\mathrm{BC}_M(Y)$. The \emph{bicycle matroid}
$\mathrm{BM}_M(Y)$ of binary matroid $M$ \emph{relative} to $Y \subseteq V$ is the (unique)
binary matroid with ground set $V$ and cycle space $\mathrm{BC}_M(Y)$. The notion of bicycle matroid for the case $Y = V$ was introduced in \cite{Jaeger/1983/SymmReprBinMat}.

In contrast to the binary case, $\mathrm{BC}_{\ker(B)}$ and $\mathrm{BC}_{\ker(B')}$ may differ when $B$ and $B'$ are different representations over $\gf{4}$ of a quaternary matroid $M$. However, we know from \cite{DBLP:journals/jct/Vertigan98} that, for all $L \subseteq \gf{4}^V$, the \emph{dimensions} of $\mathrm{BC}_{\ker(B)}$ and $\mathrm{BC}_{\ker(B')}$ are equal. We now extend this result by showing that the \emph{matroids} of $\mathrm{BC}_{\ker(B)}$ and $\mathrm{BC}_{\ker(B')}$ are equal. In fact, we show $M(\mathrm{BC}_{\ker(B)}(Y)) = M(\mathrm{BC}_{\ker(B')}(Y))$ for all $Y \subseteq V$. As a consequence we may speak of the bicycle matroid of a quaternary matroid $M$ relative to $Y \subseteq V$.
Moreover we give an explicit formula for this bicycle matroid in terms of $M$ (independent of representation). Also, the proof of this result below is direct, and therefore not obtained as a consequence of an evaluation of the Tutte polynomial as in \cite{DBLP:journals/jct/Vertigan98}.

First we prove a technical lemma.
\begin{lemma} \label{lem:matrix_bicycle}
Let $B$ be a $X \times V$-matrix over some $\field$ with characteristic $2$ such that $X \subseteq V$ and $B[X]$ the identity matrix of suitable size. Let $\alpha$ be an involutive automorphism of $\field$ and let $A = R(B,\alpha)$. Then for all $Y \subseteq V$, the null space of the $\alpha$-symmetric matrix $A+(X \cup Y)*(X\setminus Y)$ is equal to $\{v \in \ker(B) \mid \pi_Y(v) \in \alpha(\ker(B)^{\perp})\}$.
\end{lemma}
\begin{Proof}
Recall from, e.g., \cite[Proposition~2.2.23]{Oxley/MatroidBook-2nd}, that the null space of
$$
B = \bordermatrix{& \scriptstyle X& \scriptstyle V\setminus X\cr
 \scriptstyle X & I & S
}
$$
is the orthogonal complement of the null space of
$$
B' = \bordermatrix{& \scriptstyle X& \scriptstyle V\setminus X\cr
 \scriptstyle V\setminus X & -S^T & I
}.
$$
Although signs are irrelevant over fields with characteristic $2$, we leave them for didactical purposes. We observe that $B = (A+V)[X,V]$ and $\alpha(B') = (A+V)[V\setminus X,V]$. Thus, $\ker(A+V) = \ker(B) \cap \ker(\alpha(B')) = \ker(B) \cap \alpha(\ker(B)^{\perp})$ (which proves the case $Y = V$). Let
$$
S =
\bordermatrix{&\scriptstyle Y \setminus X&\scriptstyle V\setminus (Y \cup X)\cr
\scriptstyle X \setminus Y & S_1 & S_2 \cr
\scriptstyle X \cap Y & S_3 & S_4
}.
$$
Then the null space of
$$
B'' = \bordermatrix{&\scriptstyle X\setminus Y&\scriptstyle X \cap Y &\scriptstyle Y \setminus X &\scriptstyle V\setminus (Y \cup X)\cr
\scriptstyle X\setminus Y & I & 0 & S_1 & S_2 \cr
\scriptstyle X \cap Y & 0 & I & S_3 & S_4 \cr
\scriptstyle Y \setminus X & 0 & \alpha(-S_3^T) & I & 0\cr
\scriptstyle V\setminus (Y \cup X) & 0 & \alpha(-S_4^T) & 0 & 0
}.
$$
is equal to $\{v \in \ker(B) \mid \pi_Y(v) \in \alpha(\ker(B)^{\perp})\}$.
We show that the null space of $B''$ is equal to the null space of $A' = A+(X \cup Y)*(X\setminus Y)$. We have
$$
A' =
\bordermatrix{&\scriptstyle X\setminus Y&\scriptstyle X \cap Y &\scriptstyle Y \setminus X &\scriptstyle V\setminus (Y \cup X)\cr
\scriptstyle X\setminus Y & I & 0 & -S_1 & -S_2 \cr
\scriptstyle X \cap Y & 0 & I & S_3 & S_4 \cr
\scriptstyle Y \setminus X & \alpha(-S_1^T) & \alpha(-S_3^T) & I+\alpha(S_1^T) S_1 & \alpha(S_1^T) S_2 \cr
\scriptstyle V\setminus (Y \cup X) & \alpha(-S_2^T) & \alpha(-S_4^T) & \alpha(S_2^T) S_1 & \alpha(S_2^T) S_2
}.
$$
Now consider the following nonsingular matrix
$$
A'' = \bordermatrix{&\scriptstyle X\setminus Y&\scriptstyle X \cap Y &\scriptstyle Y \setminus X &\scriptstyle V\setminus (Y \cup X)\cr
\scriptstyle X\setminus Y & I & 0 & 0 & 0 \cr
\scriptstyle X \cap Y & 0 & I & 0 & 0 \cr
\scriptstyle Y \setminus X & \alpha(S_1^T) & 0 & I & 0\cr
\scriptstyle V\setminus (Y \cup X) & \alpha(S_2^T) & 0 & 0 & I
}.
$$
Observe that $A''A' = B''$ (here we use that $\field$ has characteristic $2$ to remove the ``incorrect'' signs of $(A''A')[X\setminus Y,V \setminus X]$). Hence $\ker(A+(X \cup Y)*(X\setminus Y)) = \ker(A') = \ker(B'')$.
\end{Proof}
We note that (the proof of) Lemma~\ref{lem:matrix_bicycle} holds also when automorphism $\alpha$ is not involutive. Of course, in that case, $A+(X \cup Y)*(X\setminus Y)$ need not be $\alpha$-symmetric.

It follows from the strong principal minor theorem \cite[Theorem~2.9]{Kodiyalam_Lam_Swan_2008} that, if $A$ is $\alpha$-symmetric matrix, then matroid $\max(\mathcal{M}_{A})$ is equal to $M(A)$ (in fact, the strong principal minor theorem holds for so-called ``quasi-symmetric matrices'', but it is straightforward to verify that  $\alpha$-symmetric matrices are quasi-symmetric). We use this observation in the next result.

\begin{theorem} \label{thm:bicyclemat_char}
Let $M$ be a quaternary matroid, and let $B$ be a representation of $M$ over $\gf{4}$. For all $Y \subseteq V$, $M(\mathrm{BC}_{\ker(B)}(Y))$ is equal to the matroid $\max(M+Y)$.
\end{theorem}
\begin{Proof}
Without loss of generality we assume that $B$ is a standard representation of $M$. Say, $B$ is a $X \times V$-matrix. Let $A = R(B,\invm)$ and let $A' = A+(X \cup Y)*(X\setminus Y)$. By Lemma~\ref{lem:matrix_bicycle}, $M(A') = M(\mathrm{BC}_{\ker(B)}(Y))$ and thus it suffices to show that $M(A') = \max(M+Y)$.

Recall that (through the strong principal minor theorem) $M(A') = \max(\mathcal{M}_{A'})$ since $A'$ is an $\invm$-symmetric matrix. By Theorem~\ref{thm:inv_symm_PU}, $A$ is PU and by Proposition~\ref{prop:pu_loopc}, $\mathcal{M}_{A'} = \mathcal{M}_{A}+(X \cup Y)*(X\setminus Y)$. By Proposition~\ref{prop:bouchet_twist_matroid}, $\mathcal{M}_A*X = M$ and thus $\mathcal{M}_{A}+(X \cup Y)*(X\setminus Y) = M*X+(X \cup Y)*(X\setminus Y)$.

Recall that $*$ and $+$ commute on distinct elements and that they form $S_3$ on common elements. Hence, $M*X+(X \cup Y)*(X\setminus Y) = M*(X \setminus Y)*(X \cap Y)+Y+(X\setminus Y)*(X\setminus Y) = M*(X \cap Y)+Y\dual (X\setminus Y) = M+Y \dual X$.
We obtain $M(A') = \max(M+Y \dual X)$. Since $\max(N\dual Z) = \max(N)$ for all set systems $N$ and all subsets $Z$ of the ground set, we have $\max(M+Y \dual X) = \max(M+Y)$ and the result follows.
\end{Proof}

Theorem~\ref{thm:bicyclemat_char} suggests the following definition. For a vf-safe matroid $M$, we call matroid $\max(M+Y)$ the \emph{bicycle matroid relative to} $Y \subseteq V$ and denote it by $\mathrm{BM}_M(Y)$.  Note that this definition is consistent with (and therefore generalizes) the above definition of bicycle matroid for binary matroids.

Note that $\max(M^*+V) = \max(M*V+V)= \max(M+V\dual V) = \max(M+V)$, so the bicycle matroid of $M$ is invariant under taking the dual matroid.

We call the nullity of $\mathrm{BM}_M(Y)$ the \emph{bicycle dimension} of $M$ relative to $Y$. In case $Y = V$, we simply speak of the bicycle dimension of $M$. We have $\nu(\max(M+Y)) = d_{M+Y*V} =
d_{M*V\dual Y}$, and $d_{M*V\dual V} = d_{M\dual V+V} = d_{M\dual V}$, and so we obtain the following corollary to Theorem~\ref{thm:bicyclemat_char}.
\begin{corollary} \label{cor:dim_bicycle_space}
Let $M$ be a quaternary matroid and $Y \subseteq V$. The bicycle dimension of $M$ relative to $Y$ is $d_{M*V\dual Y}$. In particular, the bicycle dimension of $M$ is
$d_{M\dual V}$.
\end{corollary}

We may rephrase Corollary~\ref{cor:dim_bicycle_space} in terms of subspaces of $\gf{4}^V$. Let $L$ be a subspace of $\gf{4}^V$. Then, by Corollary~\ref{cor:dim_bicycle_space}, $\dim(\mathrm{BC}_L(Y)) = d_{M(L)*V\dual Y}$. In particular, $\dim(L \cap \invm(L^{\perp})) = d_{M(L)\dual V}$.

The equality of the bicycle dimension and the
value $d_{M\dual V}$ was already shown for the case of binary matroids $M$ in
\cite{BH/InterlacePolyDM/14} as a consequence of calculating
the Tutte polynomial at $(-1,-1)$ in an alternative way. The
present paper explains this equality for quaternary matroids in general in a direct way (without
considering the Tutte polynomial) as a corollary to
Theorem~\ref{thm:bicyclemat_char}.

\begin{Remark}
Theorem~\ref{thm:bicyclemat_char} identifies, for quaternary matroids $M$, a relationship between the matroids $M = \max(M)$, $M^* = M*V = \max(M*V)$, and $\max(M+V)$. It turns out that the matroids $\max(M)$, $\max(M*V)$, and $\max(M+V)$ are also in some weak sense related for \dmatroids $M$ in general. Matroids $M_1$ and $M_2$ are said to be \emph{orthogonal} if for all circuits $C_1$ of $M_1$ and $C_2$ of $M_2$, $|C_1 \cap C_2| \neq 1$. It is well known that any matroid $M$ is orthogonal to its dual $M^*$. In \cite{BH/InterlacePolyDM/14}, vf-safe \dmatroids are shown to be ``essentially'' equivalent to a particular class of multimatroids \cite{DBLP:journals/siamdm/Bouchet97} called tight $3$-matroids \cite{DBLP:journals/ejc/Bouchet01} (we will not recall multimatroids in this paper). Theorem~3.2 of \cite{DBLP:journals/combinatorics/Bouchet98} shows that the matroids corresponding to disjoint transversals of a multimatroid are orthogonal when projecting the ground sets onto a common ground set $V$. This translates to vf-safe \dmatroids as follows: for any vf-safe \dmatroid $M$, the matroids $\max(M)$, $\max(M*V)$, and $\max(M+V)$ are mutually orthogonal. In case $M$ is a vf-safe matroid, we have that $M$, its dual $M^*$, and $\max(M+V)$ are mutually orthogonal.
\end{Remark}

By definition, the \dmatroid $M\dual V$ is constructed from $M$
by adding the sets that are included in an odd number of bases.
Corollary~\ref{cor:dim_bicycle_space} opens the possibility of
parity-type characterizations of the bicycle dimension. Indeed,
quaternary matroid $M$ has an odd number of bases iff the bicycle dimension of $M$ is zero, a result shown by
Chen~\cite{Chen/VectorSpaceGraph/1971} for the case where $M$
is a graphic matroid (and later realized to hold for binary matroids in general). Moreover, by the definition of dual
pivot, $d_{M\dual V} > 1$ iff the number of bases of $M$ is even
and for all $v \in V$, $v$ is in an even number of bases of $M$
iff for all $v \in V$, $v$ is in an even number of bases and in
an even number of cobases of $M$, the latter of which is the
$q>1$ characterization ($q$ being equal to the bicycle
dimension of $M$) of de~Fraysseix
\cite[Th\'eor\`eme~1]{Fraysseix/BicycleDim/78}.

\begin{Remark}
Unfortunately, the other two characterizations for $q=0$ and
$q=1$ stated in \cite[Th\'eor\`eme~1]{Fraysseix/BicycleDim/78}
do not hold. These characterizations are formulated in terms of
the principal tripartition which we recall in
Subsection~\ref{subs:tripartition}. We give a counterexample
for each characterization. The cycle matroid $M$ of $K_4$, the
complete graph on four vertices, has 16 bases, while every
element of the ground set occurs in 8 bases. Every element is
part of a $4$-cycle, which is a cocycle as well. Hence, the
tripartition equals $(P,Q,R) = (\emptyset,\emptyset,V)$. The
first characterization of
\cite[Th\'eor\`eme~1]{Fraysseix/BicycleDim/78} predicts $q=0$,
while actually $q=2$. The uniform
matroid $U_{2,3}$ is the cycle matroid of $K_3$. It has three
bases, and each element in the ground set occurs in two bases.
Moreover the ground set forms a cycle that becomes a cocycle
when any element is removed. Thus the tripartition for
$U_{2,3}$ equals $(P,Q,R) = (V,\emptyset,\emptyset)$. Now $R$
is empty, as is the set of elements occurring in an odd number
of bases. The second characterization of
\cite[Th\'eor\`eme~1]{Fraysseix/BicycleDim/78} predicts $q=1$,
while actually $q=0$. In fact, the two
characterizations for $q=0$ and $q=1$ are not disjoint, as in
the latter example also $Q=\emptyset$, predicting also $q=0$.
\end{Remark}

\section{Consequences}

In this section we discuss a number of consequences for binary
matroids of Theorem~\ref{thm:bicyclemat_char} and we give an
example. In Section~\ref{sect:penrose} we use the result to
generalize the Penrose polynomial to \dmatroids.

\subsection{Fundamental graph of a matroid}
We now consider fundamental graphs of binary matroids.
\begin{corollary}
Let $G$ be a fundamental graph of a binary matroid $M$. Then
the column matroid of $G+V$ is equal to the bicycle matroid of
$M$.
\end{corollary}
\begin{Proof}
Let $\mathcal{M}_G = M*Z$ for some $Z \in M$. Then by
Theorem~\ref{thm:bicyclemat_char}, the bicycle matroid of $M$
relative to $V$ is $\max(M+V) = \max(\mathcal{M}_G*Z+V) =
\max(\mathcal{M}_G+V\dual Z) = \max(\mathcal{M}_G+V) =
\max(\mathcal{M}_{G+V})$ which in turn is equal to the column
matroid of $G+V$.
\end{Proof}
Consequently, if $G_1$ and $G_2$ are fundamental graphs of some
binary matroid $M$, then the column matroids of $G_1+V$ and
$G_2+V$ are equal.

Since every bipartite graph is the fundamental graph of some
matroid $M$, we obtain the following result stated (without
proof) in \cite{Jaeger/1983/SymmReprBinMat}.
\begin{corollary}[Proposition~3 of \cite{Jaeger/1983/SymmReprBinMat}]
Let $M$ be binary matroid. Then $M$ is the column matroid of a
graph $G$ such that $G+V$ is bipartite iff $M$ is the bicycle
matroid of a binary matroid.
\end{corollary}

\begin{Remark}
The Tanner graph \cite{DBLP:journals/tit/Tanner81} is a popular
notion within coding theory. A (linear) \emph{code} $\mathcal{C}$ is a
subspace of $\gf2^V$ (for some finite set $V$), a
\emph{parity-check} matrix $H$ for $\mathcal{C}$ is a matrix
with $\ker(H) = \mathcal{C}$. Matrix $H$ is said to be in a
\emph{standard form} if $H = (I \quad B)$ where $I$ is an
identity matrix. If $H$ is a $m \times n$-parity-check matrix
in standard form, then the \emph{Tanner graph} $T$ of $H$ is a
$(U,V)$-bipartite graph with $U=\{u_1,\ldots,u_m\}$ and
$V=\{v_1,\ldots,v_n\}$ such that $\{u_i,v_j\}$ is an edge of
$T$ iff entry $H_{i,j}$ is equal to $1$. The elements of $U$
and $V$ are called \emph{check nodes} and \emph{bit nodes},
respectively. The vertices $v_1,\ldots,v_m$ and their edges are
often ignored, see, e.g., Fig.~2 in
\cite{DBLP:conf/isit/KnudsenRDPR09}. The obtained bipartite
graph $G$ is therefore exactly a fundamental graph of the
binary matroid $M$ with cycle space $\mathcal{C}$. Hence it
seems fruitful to consider Tanner graphs from a matroid
perspective (indeed, e.g., edge local complementation is
applied to the Tanner graph, see
\cite{DBLP:conf/isit/KnudsenRDPR09}, which corresponds to
taking a different fundamental graph of $M$). However,
surprisingly, considering the extensive literature on the
notion of Tanner graph, this perspective seems to have not yet
been taken.
\end{Remark}

In the next subsections we meet the fundamental graph again.

\subsection{Principal tripartition} \label{subs:tripartition}
The well-known principal tripartition for binary
matroids from Rosenstiehl and Read \cite[Theorem~2.1]{Rosenstiehl1978195} is as follows.
\begin{proposition}[Principal tripartition \cite{Rosenstiehl1978195}] \label{prop:tripartition}
Let $M$ be a binary matroid. Then every $u \in V$ belongs to exactly one of the following sets: $P = \{v \in V \mid \exists X \in \mathcal{CS}_M, v \in X, X \setminus \{v\} \in \mathcal{CS}_M^{\perp}\}$, $Q = \{v \in V \mid \exists X \in \mathcal{CS}_M^{\perp}, v \in X, X \setminus \{v\} \in \mathcal{CS}_M\}$, and $R = \{v \in V \mid \exists X \in \mathcal{CS}_M \cap \mathcal{CS}_M^{\perp}, v \in X \}$
\end{proposition}

We now generalize the principal tripartition result from binary matroids to vf-safe matroids.
The characterization of the previous section allow us to use a formulation that avoids (bi)cycles.
\begin{theorem} \label{thm:charPQR}
Let $M$ be a vf-safe \dmatroid. Then every element of $V$ belongs to exactly one of the following sets: $P = \{v \in V \mid v \text{ is not a
coloop of } \max(M+V+v) \}$, $Q = \{v \in V \mid v \text{ is not a
coloop of } \max(M+V*v) \}$, and $R = \{v \in V \mid v \text{ is not a
coloop}\allowbreak\text{of } \max(M+V) \}$. Moreover, this tripartition coincides with the tripartition of Proposition~\ref{prop:tripartition} when $M$ is a binary matroid.
\end{theorem}
Note that, in particular, quaternary matroids satisfy the tripartition result of Theorem~\ref{thm:charPQR}.

The proof of Theorem~\ref{thm:charPQR} uses the following result from \cite{BH/NullityLoopcGraphsDM/10}.
\begin{proposition}[Theorem~14 in
\cite{BH/NullityLoopcGraphsDM/10}] \label{prop:mmmplusone} Let
$M$ be a vf-safe \dmatroid and let $v \in V$. Then the matroids $\max(M)$, $\max(M*v)$, and $\max(M+v)$ are such that precisely two of the three are equal, to say $M_1$. Moreover, the rank of the third $M_2$ is one smaller than the rank of $M_1$, and $M_1$ is the direct sum of $M_2 \vertexrem v$ and the matroid consisting of the coloop $v$.
\end{proposition}

\begin{Proof} \emph{(of Theorem~\ref{thm:charPQR})}
By applying Proposition~\ref{prop:mmmplusone} to $M+V$, we have that two of three matroids $\max(M+V+v)$, $\max(M+V*v)$, $\max(M+V)$ are equal and have $v$ as a coloop and the third does not have $v$ as a coloop. Hence $v$ belongs to precisely one of $P$, $Q$, and $R$.

Assume now that $M$ is a binary matroid. We show that the tripartitions of Proposition~\ref{prop:tripartition} and Theorem~\ref{thm:charPQR} coincide. Now, there exists a $X \in \mathcal{CS}_M$ with $X \setminus \{v\} \in \mathcal{CS}_M^{\perp}$ and $v \in X$ iff there exists a $X \in \mathcal{CS}_{\mathrm{BM}_M(V\setminus \{v\})}$ with $v \in X$ iff $v$ is not a coloop of $\mathrm{BM}_M(V\setminus \{v\}) = \max(M+(V\setminus \{v\})) = \max(M+V+v)$ (where the first equality is by Theorem~\ref{thm:bicyclemat_char}).
Similarly, there exists a $X \in \mathcal{CS}_M^{\perp} = \mathcal{CS}_{M^*}$ with $X \setminus \{v\}  \in \mathcal{CS}_M = \mathcal{CS}_{M^*}^{\perp}$ and $v \in X$ iff there exists a $X \in \mathcal{CS}_{\mathrm{BM}_{M^*}(V\setminus \{v\})}$ with $v \in X$ iff $v$ is not a coloop of $\mathrm{BM}_{M^*}(V\setminus \{v\}) = \max(M*V+(V\setminus \{v\})) = \max(M+(V\setminus \{v\}) *v \dual (V\setminus \{v\})) = \max(M+(V\setminus \{v\}) *v \dual v) = \max(M+V*v)$.
Finally, there exists a $X \in \mathcal{CS}_M \cap \mathcal{CS}_M^{\perp} = \mathcal{CS}_{\mathrm{BM}_M(V)}$ with $v \in X$ iff $v$ is not a coloop of $\mathrm{BM}_{M}(V) = \max(M+V)$.
\end{Proof}

Theorem~\ref{thm:charPQR} and its relation to Proposition~\ref{prop:mmmplusone} allows one to generalize results associated to the principal tripartition result (and, moreover, allows for easier proofs of these results). For example Table~3 in \cite{Rosenstiehl1978195}, which states how the tripartition changes for a cycle matroid of a graph $G$ when applying various operations on $G$, is readily obtained as a consequence of Proposition~\ref{prop:mmmplusone}.

In \cite{BT/AdjMatroidGraph/2011}, the graph counterpart of
Proposition~\ref{prop:mmmplusone} (i.e., the case $M =
\mathcal{M}_G$ for some graph $G$) was explicitly seen as a
tripartition result with the property of ``being a coloop''
that decides to which of the three classes of the tripartition
a particular vertex belongs (similar as in
Theorem~\ref{thm:charPQR}). However, no concrete link with the
result of \cite{Rosenstiehl1978195} was established in
\cite{BT/AdjMatroidGraph/2011}.

Let, for a vf-safe \dmatroid $M$ and $v \in V$, $\mathrm{nmax}_M(v) = \max\{\nu(\max(M)),\allowbreak\nu(\max\allowbreak(M*v)),\allowbreak\nu(\max(M+v))\}$. Note that by using Proposition~\ref{prop:mmmplusone}, the definitions of $P$, $Q$, and $R$ of Theorem~\ref{thm:charPQR} can be rephrased as follows: $v \in P$ iff $\nu(\max(M+V+v)) = \mathrm{nmax}_{M+V}(v)$, $v \in Q$ iff $\nu(\max(M+V*v)) = \mathrm{nmax}_{M+V}(v)$, and $v \in R$ iff $\nu(\max(M+V)) = \mathrm{nmax}_{M+V}(v)$. One may again rephrase this in terms of rank instead of nullity. However, we choose nullity due to the following result.

We let $\nu(G)$ be the nullity of the adjacency matrix of a graph $G$. Moreover, let for all $v \in V$,
$\mathrm{nmax}_G(v) = \max\{\;\nu(G), \,\nu(G\vertexrem v), \,\nu(G+v)\;\}$. We reformulate the tripartition for binary matroids $M$ in terms of a fundamental graph of $M$.
\begin{corollary} \label{cor:tripartition_fund_graph}
Let $P,Q,R$ be the tripartition associated with a binary matroid $M$. Let $G$ be a fundamental graph of $M = \mathcal{M}_G*Z$. If $v \in V \setminus Z$, then
\begin{enumerate}
\item $v \in P$ iff $\mathrm{nmax}_{G+V}(v) = \nu(G+V+v))$,
\item $v \in Q$ iff $\mathrm{nmax}_{G+V}(v) = \nu(G+V\vertexrem v)$, and
\item $v \in R$ iff $\mathrm{nmax}_{G+V}(v) = \nu(G+V)$.
\end{enumerate}
If $v \in Z$, then the roles of $P$ and $Q$ are reversed.
\end{corollary}
\begin{Proof}
Let $v \in V \setminus Z$. Then $\max(M+V+v) =
\max(\mathcal{M}_G*Z+V+v) = \max(\mathcal{M}_G+V+v) =
\max(\mathcal{M}_{G+(V \setminus \{v\})})$. Hence $\nu(\max(M+V+v)) =
\nu(G+(V \setminus \{v\}))$. Similarly, $\nu(\max(M+V*v)) =
\nu(G+V*v)$ and $\nu(\max(M+V)) = \nu(G+V)$. It is a
well-known property of the Schur complement
$F*v\vertexrem v$ for a matrix $F$ that its nullity is equal to the nullity of $F$,
see, e.g., \cite{SchurBook2005}. Hence $\nu(G+V*v) =
\nu(G+V\vertexrem v)$.

Finally, let $v \in Z$. Then $\max(M+V+v) =
\max(\mathcal{M}_G*Z+V+v) =
\max(\mathcal{M}_G+V\dual v+v) =
\max(\mathcal{M}_G+V*v)$ and similarly $\max(M+V*v) =
\max(\mathcal{M}_G+V+v)$. Hence the roles of $P$ and $Q$
are reversed with respect to the case of $v \in V \setminus Z$.
\end{Proof}

The fact that the three values $\nu(G+v)$, $\nu(G\vertexrem v)$, $\nu(G)$ in Corollary~\ref{cor:tripartition_fund_graph} are of the
form $m,m,m+1$ (in some order)
has been shown in \cite[Lemma~2]{DBLP:journals/ejc/BalisterBCP02}, see also
\cite{BH/PivotLoopCompl/09}. As a consequence, it suffices to know two of
these three values to be able to determine the third.

\subsection{Eulerian and bipartite binary matroids}

Matroid $M$ is said to be \emph{bipartite} when every circuit
of $M$ is of even cardinality, and $M$ is said to be
\emph{Eulerian} when there are disjoint circuits of $M$ whose union is equal to $V$. If $M$ is binary, then $M$ is Eulerian iff $V \in \mathcal{CS}_M$. It is shown in
\cite{Welsh/EulerBipartiteDual} that a binary matroid $M$ is
Eulerian iff its dual $M^*$ is bipartite. We obtain the
following two dual characterizations.

\begin{theorem} \label{thm:bipartite_char}
Let $M$ be a binary matroid. \\
{\rm(1)}
$M$ is bipartite iff $M+V$ is an even \dmatroid.
\\
{\rm(2)}
$M$ is Eulerian iff $M\dual
V$ is an even \dmatroid.
\end{theorem}
\begin{Proof}
(1)
By the proof of
\cite[Proposition~2]{Jaeger/1983/SymmReprBinMat}, $M$ is
bipartite iff each diagonal entry of $G+X*X$ is
$1$ where $G$ is the fundamental graph of $M$ with respect to some $X
\in M$. Since $\mathcal{M}_{G} = M*X$, $M$ is bipartite iff $M
\dual X +V = M+V*X$ is an even \dmatroid. The latter is in turn
equivalent to $M+V$ being an even \dmatroid.

(2)
Note $M^*+V=M*V+V=M\dual V*V$ is even iff $M\dual V$ is even.
\end{Proof}

\begin{Remark}
Theorem~\ref{thm:bipartite_char} suggests an alternative extension of the
notions of bipartite and Eulerian from binary matroids to
vf-safe (delta-)matroids. This alternative extension is natural as the two
notions remain each others dual notions---the (original)
notions of bipartite and Eulerian are known to not be dual for
(nonbinary) matroids in general. For example, the vf-safe (in fact, quaternary) uniform matroid $U_{3,6}$ is bipartite, but $U^*_{3,6} = U_{3,6}$ is not Eulerian.
\end{Remark}

\begin{Remark}
Proposition~1 of \cite{Jaeger/1983/SymmReprBinMat} (which
is used in the proof of
\cite[Proposition~2]{Jaeger/1983/SymmReprBinMat}) states that a
matroid $M$ is binary iff $M = \max(\mathcal{M}_{F})$ for some
graph $F$. In fact, the proof of this result implicitly takes
$F = G+X*X$ where $G$ is the fundamental graph of $M$ with respect to
$X$. Indeed, $M = \max(M) = \max(M\dual X) =
\max(\mathcal{M}_{G+X*X})$ where $G$ is such that
$\mathcal{M}_G = M*X$.
\end{Remark}


\subsection{Example}
\def\VR{\kern-\arraycolsep\strut\vrule &\kern-\arraycolsep}
\def\vr{\kern-\arraycolsep & \kern-\arraycolsep}

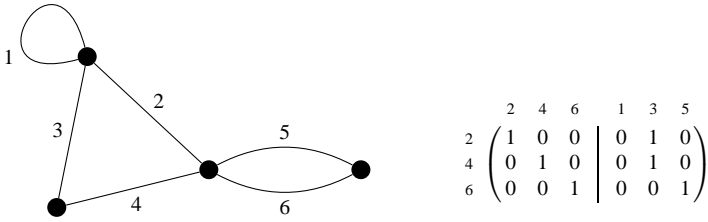
\begin{figure}
\begin{center}
\begin{picture}(50,35)(-7,-5)
\gasset{AHnb=0,Nw=2.5,Nh=2.5,Nframe=n,Nfill=y}
  \node(1)(04,20){}
  \node(2)(00,00){}
  \node(3)(20,05){}
  \node(4)(40,05){}
  \drawloop[loopangle=150,ELpos=30](1){1}
  \drawedge(1,3){2}
  \drawedge(2,1){3}
  \drawedge(3,2){4}
  \drawedge[curvedepth=3](3,4){5}
  \drawedge[curvedepth=3](4,3){6}
\end{picture}%
\hspace{1cm}
\raisebox{10mm}{\(
\bordermatrix{
   & \scriptstyle 2 & \scriptstyle 4 & \scriptstyle 6 & \vr \scriptstyle 1 & \scriptstyle 3 & \scriptstyle 5 \cr
\scriptstyle 2  & 1 & 0 & 0 & \VR 0 & 1 & 0 \cr
\scriptstyle 4  & 0 & 1 & 0 & \VR 0 & 1 & 0 \cr
\scriptstyle 6  & 0 & 0 & 1 & \VR 0 & 0 & 1
}
\)}
\end{center}
\caption{A graph $F$ and a binary representation of its cycle matroid $M$.}
\label{fig:ex-graph-matroid}
\end{figure}

Let us consider the graph $F$ of Figure~\ref{fig:ex-graph-matroid} (left-hand side) with six (labeled) edges, and consider the cycle matroid $M$ of $F$ over $V = \{ 1, \dots,
6\}$, see Figure~\ref{fig:ex-graph-matroid} (right-hand side). To avoid notational clutter, we often denote sets within sets by juxtaposition in this example. The bases of $M$ are the six spanning trees of $F$, thus $M = ( V, \{ 235, 236, 245, 246, 345, 346 \} )$.

The cycle space of $M$ has dimension $3$, and is generated by
$\{ 1, 234, 56 \}$. Its cocycle space is also of dimension $3$,
and is generated by $\{ 23, 24, 56 \}$.

The empty set $\emptyset$ is not a set in $M\dual V$, as $M$
has an even number of bases. However, $M\dual V$ contains
$\{5\}$ and $\{6\}$, as both are contained in three bases.
Thus the bicycle dimension of $M$ is $d_{M \dual V} = 1$.

From $M$ one constructs $M+V$ by adding the sets that contain
an odd number of bases, these are eight $4$-element sets
and the $5$-element
sets $\{ 1,2,3,4,5\}$ and $\{1,2,3,4,6\}$. Thus matroid $\max(M+V)$ equals
$(V, \{ 12345,\allowbreak 12346  \})$. Hence, the only nontrivial cycle
of $\max(M+V)$ is $\{ 5,6 \}$.

Consider basis $Z = \{2,4,6\}$ of $M$. The standard
representation of $M$ with respect to $Z$ is given in
Figure~\ref{fig:ex-graph-matroid}.
Using that representation we deduce the (bipartite) fundamental
graph $G$ of $M$ with respect to $Z$ (the construction is described in
\cite{bouchet1987}).

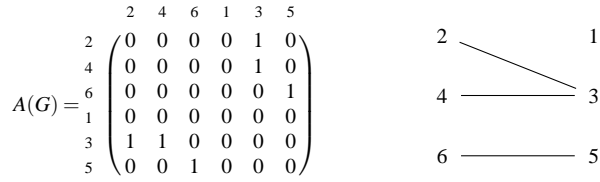
\begin{figure}
\begin{center}
\( A(G) =
\bordermatrix{
   & \scriptstyle 2 & \scriptstyle 4 & \scriptstyle 6 & \scriptstyle 1 & \scriptstyle 3 & \scriptstyle 5 \cr
\scriptstyle 2  & 0 & 0 & 0 & 0 & 1 & 0 \cr
\scriptstyle 4  & 0 & 0 & 0 & 0 & 1 & 0 \cr
\scriptstyle 6  & 0 & 0 & 0 & 0 & 0 & 1 \cr
\scriptstyle 1  & 0 & 0 & 0 & 0 & 0 & 0 \cr
\scriptstyle 3  & 1 & 1 & 0 & 0 & 0 & 0 \cr
\scriptstyle 5  & 0 & 0 & 1 & 0 & 0 & 0
}
\)
\hspace{1cm}
\raisebox{-8mm}{%
\begin{picture}(30,20)(-5,0)
\gasset{AHnb=0,Nw=5,Nh=5,Nframe=n,Nfill=n}
  \node(2)(00,18){2}
  \node(4)(00,10){4}
  \node(6)(00,02){6}
  \node(1)(20,18){1}
  \node(3)(20,10){3}
  \node(5)(20,02){5}
  \drawedge(2,3){}
  \drawedge(4,3){}
  \drawedge(6,5){}
\end{picture}%
}
\end{center}
\caption{The fundamental graph of $M$ with respect to basis $\{2,4,6\}$.}
\label{fig:ex-fundamental}
\end{figure}

The \dmatroid ${\cal M}_G$ is equal to $M*Z = ( V, \{ 3456, 34,
56, \emptyset, 2356, 23 \} )$.

The adjacency  matrix of $G+V$ is as follows. Due to the simple
block structure of the matrix, nullities are easily computed.
From that we infer the tripartition using
Corollary~\ref{cor:tripartition_fund_graph}.

\noindent
\begin{center}
$A(G+V) =
\bordermatrix{
   & \scriptstyle 1 & \scriptstyle 2& \scriptstyle 3 & \scriptstyle 4 & \scriptstyle 5 & \scriptstyle 6 \cr
\scriptstyle 1  & 1 & 0 & 0 & 0 & 0 & 0 \cr
\scriptstyle 2  & 0 & 1 & 1 & 0 & 0 & 0 \cr
\scriptstyle 3  & 0 & 1 & 1 & 1 & 0 & 0 \cr
\scriptstyle 4  & 0 & 0 & 1 & 1 & 0 & 0 \cr
\scriptstyle 5  & 0 & 0 & 0 & 0 & 1 & 1 \cr
\scriptstyle 6  & 0 & 0 & 0 & 0 & 1 & 1
}
$
\hspace{5mm}
$
\begin{array}{r|*6c}
      v                & \;1\; & \;2\; & \;3\; & \;4\; & \;5\; & \;6\; \\\hline
\nu(G+V+v)         & 2 & 1 & 2 & 1 & 0 & 0 \\
\nu(G+V\setminus \{v\}) & 1 & 2 & 1 & 2 & 0 & 0 \\
\nu(G+V)               & 1 & 1 & 1 & 1 & 1 & 1 \\
\mbox{in } Z           & \mbox{no} & \mbox{yes} & \mbox{no} & \mbox{yes} & \mbox{no} & \mbox{yes} \\[1ex]
\mbox{tripartition}    & P & P & P & P & R & R
\end{array}
$
\end{center}

\section{Penrose Polynomial} \label{sect:penrose}
Jaeger \cite{Jaeger1990Transition} defines the \emph{Penrose
polynomial} for $4$-regular graphs, inspired by the work of
Penrose \cite{Penrose/Tensors/1971}, and Aigner and Mielke
\cite{AignerPenroseBinMat} show that this notion can be defined more
generally for a binary matroid $M$ as
$$
P_M(y) = \sum_{X \subseteq V} (-1)^{|X|} y^{\dim(\mathrm{B}_M(X))}.
$$
By Corollary~\ref{cor:dim_bicycle_space} we obtain
$$
P_M(y) = \sum_{X \subseteq V} (-1)^{|X|}
y^{d_{M*V\dual X}}.
$$
The latter formulation allows us to consider $P_M(y)$ for arbitrary vf-safe \dmatroids $M$ (or, indeed, set systems) instead of binary matroids.
We first consider the transition polynomial for set systems,
which has the generalized Penrose polynomial as a
specialization.

In this section all set systems are assumed to be proper.
Note that $M_\emptyset=(\emptyset,\{\emptyset\})$ is proper.

\subsection{Transition Polynomials}
In a $4$-regular graph we can partition the edges into a set of
circuits (called a circuit partition). The number of resulting
circuits depends on the choices (or transitions) made at each
vertex. Jaeger shows how several classic graph polynomials
arise by counting circuits and applying particular weights for
the three possible transitions at each vertex
\cite{Jaeger1990Transition}. Here we follow this approach, but
in the more abstract fashion given in \cite{BH/InterlacePolyDM/14}.

Let $V$ be a finite set. We define $\mathcal{P}_3(V)$ to be the
set of triples $(V_1,V_2,V_3)$ where $V_1$, $V_2$, and $V_3$
are pairwise disjoint subsets of $V$ such that $V_1 \cup V_2
\cup V_3 = V$. Therefore $V_1$, $V_2$, and $V_3$ form an
``ordered partition'' of $V$ where $V_i = \emptyset$ for some
$i \in \{1,2,3\}$ is allowed.

We now recall the transition polynomial for set systems from \cite{BH/InterlacePolyDM/14}.
\begin{definition} \label{def:trans_pol}
Let $M$ be a proper set system. 
We define the \emph{transition polynomial} of $M$ (weighted by
$[a,b,c]$) as follows:
\begin{eqnarray*}
Q_{[a,b,c]}(M)(y) = \sum_{(A,B,C) \in \mathcal{P}_3(V)} a^{|A|} b^{|B|} c^{|C|}
y^{d_{M*B\dual C}}.
\end{eqnarray*}
\end{definition}

The next lemma show that $P_M(y)$ is a specialization of $Q_{[a,b,c]}(M)(y)$.
\begin{lemma}
Let $M$ be a proper set system. Then $P_M(y) = Q_{[0,1,-1]}(M)(y)$.
\end{lemma}
\begin{Proof}
Indeed, we have $Q_{[0,1,-1]}(M)(y) = \sum_{(\emptyset,B,C) \in \mathcal{P}_3(V)} (-1)^{|C|}
y^{d_{M*B\dual C}} = \sum_{C \subseteq V} (-1)^{|C|} \allowbreak y^{d_{M*(V \setminus C)\dual C}}$ and $d_{M*(V \setminus C)\dual C} = d_{M*(V \setminus C)\dual C+C}$, where in the last equality we use that $d_{N+X} = d_N$ for all set systems $N$ over $V$ and $X \subseteq V$. Now, $d_{M*(V \setminus C)\dual C+C} = d_{M*(V \setminus C)*C\dual C} \allowbreak = d_{M*V\dual C}$ and the result follows.
\end{Proof}

Several other specializations of the transition polynomial are
well-known polynomials. It is shown in
\cite{BH/InterlacePolyDM/14} that the two-variable interlace
polynomial \cite{ArratiaBS04} of a graph $G$ is equal to
$q(G)(x,y) = Q_{[1,x-1,0]}(\mathcal{M}_G)((y-1)/(x-1))$ (the
single-variable case is the case $x=2$, see
\cite{Arratia2004199}). Moreover,
$Q_{[0,b,c]}(\mathcal{M}_G)(y)$ equals the bracket polynomial
of $G$ as studied in \cite{Traldi/Bracket1/09}. Furthermore, we recall from \cite{BH/InterlacePolyDM/14} that the Tutte polynomial
$t_M(x,y)$ for a matroid $M$ is closely related to $Q_{[a,b,0]}(M)(y)$.

\begin{proposition} [Theorem~24 of \cite{BH/InterlacePolyDM/14}]\label{prop:transtition-tutte}
Let $M$ be a matroid. Then $Q_{[a,b,0]}(M)(y) = a^{\nu(M)}
b^{\rho(M)} t_M(1+\frac ab y, 1+\frac ba y)$, where $t_M$ is
the Tutte polynomial.
\end{proposition}

We quote several results from \cite{BH/InterlacePolyDM/14} that will
be useful in the present paper.
The following result
illustrates the effect on the transition polynomials $Q_{[a,b,c]}(M)$
of the application of the operations ${}+V$, ${}\dual V$, and ${}*V$.
\begin{proposition}[\cite{BH/InterlacePolyDM/14}] \label{prop:multivar_int_pol}
Let $M$ be a proper set system over $V$. Then $Q_{[a,b,c]}(M)(y) =
Q_{[a,c,b]}(M+V)(y) = Q_{[c,b,a]}(M\dual V)(y) = Q_{[b,a,c]}(M*V)(y)$.
\end{proposition}

We say that $u \in V$ is \emph{singular} in a proper set system $M$ if $u$ is a loop or a coloop of $M$.

The transition polynomial satisfies a recursive formulation \cite[Theorem~28]{BH/InterlacePolyDM/14}. Here we state the case where the $c = 0$.
\begin{proposition}[\cite{BH/InterlacePolyDM/14}]
\label{prop:recursive_rel} Let $M$ be a \dmatroid, and let
$Q(M)(y) = Q_{[a,b,0]}(M)(y)$.

\noindent {\rm(0)} For $M = M_\emptyset =
(\emptyset,\{\emptyset\})$ we have $Q(M)(y) = 1$.

Now let $u \in V$.

\noindent {\rm(1)}  If $u$ is nonsingular in $M$, then
$ Q(M)(y) = a\, Q(M\vertexrem u)(y) + b\, Q(M*u\vertexrem u)(y)$.

\noindent {\rm(2)} If $u$ is a loop of $M$, then $Q(M)(y) = (a + b
y)\,Q(M\vertexrem u)(y).$

\noindent {\rm(3)} If $u$ is a coloop of $M$, then
$Q(M)(y) = (b + a y)\,Q(M*u\vertexrem u)(y).$
\end{proposition}
Note that the recursive relation of Proposition~\ref{prop:recursive_rel} \emph{characterizes} $Q_{[a,b,0]}(M)(y)$.

Let us write for graphs $G$, $Q_{[a,b,c]}(G)(y) = Q_{[a,b,c]}(\mathcal{M}_G)(y)$.
\begin{Remark}
Theorem~\ref{thm:inv_repr_gf4_vfclosed} directly implies that the graph polynomials $Q_{[a,b,c]}(G)$ can be straightforwardly defined for $\invm$-symmetric $V \times V$-matrices $A$ over $\gf{4}$ in general instead of graphs $G$ (i.e., symmetric $V \times V$-matrices over $\gf{2}$), while maintaining their recursive relations and some of its evaluations. As an example, consider the interlace polynomial $Q(A)(y) = \sum_{X\subseteq V} \sum_{Z\subseteq X} (y-2)^{n(A+Z[X])}$ of \cite{Aigner200411} for $\invm$-symmetric $V \times V$-matrices $A$. Following the exact same reasoning as done in \cite{BH/InterlacePolyDM/14} but for $\invm$-symmetric $V \times V$-matrices $A$ instead of graphs $G$, we have that $Q(A)(y) = Q_{[1,1,1]}(\mathcal{M}_A)(y-2)$. By Theorem~\ref{thm:inv_repr_gf4_vfclosed}, $\mathcal{M}_A$ is vf-safe. Again following the reasoning of \cite{BH/InterlacePolyDM/14} by using the recursive relation of $Q_{[a,b,c]}(M)(y)$ (of which Proposition~\ref{prop:recursive_rel} above is a special case), we have that (1) if $\{u,v\} \subseteq V$ with $u \neq v$, $A[\{u\}] = A[\{v\}] = 0$, and $A[\{u\},\{v\}] \neq 0$, then $Q(A)(y) = Q(A\vertexrem u)(y) + Q(A\dual \{u\}\vertexrem u)(y) + Q(A*\{u,v\}\vertexrem u)(y)$, and (2) if $A[V,\{u\}]$ is a zero vector, then $Q(A)(y) = y Q(A\vertexrem u)(y)$. Also, we have the evaluation $Q(A)(0) = 0$ if $|V| > 0$.
\end{Remark}

We now consider the case where $M$ is even.
\begin{lemma}  \label{lem:even-p1-q1}
Let $M$ be a proper even set system. Then $Q_{[a,b,0]}(M)(y) = (-1)^{d_M} Q_{[a,-b,0]}\allowbreak(M)\allowbreak(-y) \allowbreak = (-1)^{d_{M*V}} Q_{[-a,b,0]}(M)(-y)$.
\end{lemma}
\begin{Proof}
We have $Q_{[a,b,0]}(M)(y) = \sum_{(A,B,\emptyset) \in \mathcal{P}_3(V)} a^{|A|} b^{|B|} y^{d_{M*B}}$. Since $M$ is even, the parity of $|B|+d_{M*B}$ does not depend on $B$. Hence, $Q_{[a,-b,0]}(M)(-y) = \sum_{(A,B,\emptyset) \in \mathcal{P}_3(V)} \allowbreak a^{|A|} (-b)^{|B|} (-y)^{d_{M*B}} = (-1)^{d_M} Q_{[a,b,0]}(M)(y)$. Similarly, since $M$ is even and $B = V\setminus A$, the parity of $|A|+d_{M*V*A} = |V|-|B|+d_{M*B}$ does not depend on $B$ and thus does not depend on $A$. Hence, $Q_{[-a,b,0]}(M)(-y) = \sum_{(A,B,\emptyset) \in \mathcal{P}_3(V)} (-a)^{|A|} (-b)^{|B|} (-y)^{d_{M*V*A}} \allowbreak = (-1)^{d_{M*V}} Q_{[a,b,0]}(M)(y)$.
\end{Proof}

\subsection{The specialization $Q_{[1,-1,0]}(M)(y)$}
Instead of directly studying the Penrose polynomial $P_M(y)$, we consider first $P_{M+V*V}(y)$ which by Proposition~\ref{prop:multivar_int_pol} is equal to $Q_{[0,1,-1]}(M+V*V)(y) = Q_{[1,0,-1]}(M+V)(y) = Q_{[1,-1,0]}(M)(y) = \allowbreak \sum_{X \subseteq V} \allowbreak (-1)^{|X|} y^{d_{M*X}}$. We denote $P_{M+V*V}(y) = Q_{[1,-1,0]}(M)(y)$ by $p_1(M)(y)$.
The advantage of $p_1(M)(y)$ over $P_M(y)$ is that $p_1(M)(y)$ is defined in terms of twist instead of the more elaborate loop complementation. As a result, it is easier to prove results for $p_1(M)(y)$ and then translate them to $P_M(y)$ instead of working directly with $P_M(y)$. Also, for some results concerning $p_1(M)(y)$, $M$ needs only to be \dmatroid, where the corresponding result for $P_M(y)$ requires $M$ to be a vf-safe \dmatroid.

Let $M$ be a proper set system. We obviously have, $p_1(M)(1) = 0$ if $|V| > 0$. Also, $p_1(M*X)(y) = (-1)^{|X|} p_1(M)(y)$ for all $X \subseteq V$. As a consequence, if $M$ is such that $M*X = M$ for some $X \subseteq V$ with $|X|$ odd, then $p_1(M)(y) = 0$. Note that, in this case, $M$ is not even. Finally note that by Proposition~\ref{prop:transtition-tutte}, $p_1(M)(y) = (-1)^{\rho(M)} t_M(1-y, 1-y)$ for the case where $M$ is a matroid.

\begin{Example}\label{ex:fano-p1}
Let $F_7$ be the Fano matroid, see \cite{Oxley/MatroidBook-2nd} for a description of this matroid.
We compute $p_1(F_7)$ by case analysis on the cardinality of $X \subseteq V$, from $0$ to $7$,
$p_1(F_7)(y) =  y^3 - 7 y^2 + 21 y^1 -(28 y^0 +7 y^2) + 35 y^1 - 21 y^2 + 7 y^3 - y^4 = -y^4 +8 y^3 -35 y^2 + 56 y - 28$.

We claim that $F_7\dual V = F_7$. Clearly the triplets (three-element subsets) of $F_7$ and $F_7\dual V$ coincide. It is easy
to verify that every smaller set belongs to an even number of
bases, and does not belong to $F_7\dual V$. Each two-element
subset belongs to five triplets, of which one is a line. Thus
it is contained in four bases. Each single-element subset
belongs to 15 triplets, including three lines. Thus it is
contained in twelve bases. The empty set is contained in all 28
bases.

As $F_7\dual V = F_7$ we know that $F_7 + V *V = F_7+V$, thus $p_1(F_7+V) = 0$.
We also conclude that $F_7+V$ consists of the 28 bases of $F_7$ and the 28 bases of $F_7^*$.
\end{Example}

By Proposition~\ref{prop:recursive_rel}, we have the following
recursive relation for $p_1(M)$.
\begin{corollary} \label{cor:rec_p1}
Let $M$ be a \dmatroid, and $u \in V$. If $u$ is nonsingular in $M$, then
\[
p_1(M)(y) = p_1(M\vertexrem u)(y) - p_1(M*u\vertexrem u)(y).
\]
If $u$ is a loop of $M$, then $p_1(M)(y) = (1-y) p_1(M\vertexrem u)(y)$, and if $u$ is a coloop of $M$, then $p_1(M)(y) = (y-1) p_1(M*u\vertexrem u)(y)$.
\end{corollary}

\begin{Remark}
It may strike the reader that the recursive relations of
Corollary~\ref{cor:rec_p1} are very similar to the recursive
relations of the characteristic polynomial $c_M(y)$ of a
matroid $M$. Indeed, for matroids $M$, the recursive relations
coincide when either (1) $u$ is nonsingular in $M$ or (2) $u$ is
a coloop of $M$. However, when $u$ is a loop of $M$, then
$c_M(y)$ differs as it is equal to $0$.
\end{Remark}

Let us write for graphs $G$, $p_1(G)(y) \allowbreak =
\allowbreak p_1(\mathcal{M}_G)(y)$. The next lemma shows that the graph polynomial $p_1(G)(y)$ may (in a way similar to
the interlace polynomial \cite{Arratia2004199}) be recursively
computed. We show in Theorem~\ref{thm:penrose_graphpol} below that $p_1(G)(y)$ computes (up to a sign) the Penrose polynomial for a binary matroid $M$ where $G$ is a particular graph depending on $M$. Again note that it is straightforward to consider more generally $p_1(A)(y)$ for $\invm$-symmetric $V \times V$ matrices $A$ instead of graphs $G$, and results such as Lemma~\ref{lem:p1_graph} can be formulated for $p_1(A)(y)$ as well. However, for convenience we choose to restrict to graphs $G$.
\begin{lemma} \label{lem:p1_graph}
Let $G$ be a graph. Then
$$
p_1(G)(y) = \sum_{X \subseteq V} (-1)^{|X|} y^{\nu(G[X])}.
$$
Moreover, $p_1(G)(y)$ satisfies the following characterizing recursive relation.
If $u$ is a looped vertex, then
$$
p_1(G)(y) = p_1(G\vertexrem u)(y) - p_1(G*u\vertexrem u)(y).
$$
If
$\{u,v\}$ is an edge where both $u$ and $v$ are not looped,
then
$$
p_1(G)(y) = p_1(G\vertexrem u)(y) + p_1(G*\{u,v\}\vertexrem u)(y).
$$
If $u$ is an isolated vertex (i.e., no edge is
adjacent to $u$) of $G$, then
$$
p_1(G)(y) = (1-y) p_1(G\vertexrem u)(y).
$$
Finally, if $G$ is the empty graph, then $p_1(G)(y) = 1$.
\end{lemma}
\begin{Proof}
We have $p_1(G)(y) = p_1(\mathcal{M}_G)(y) = \sum_{X \subseteq
V} (-1)^{|X|} y^{d_{\mathcal{M}_{G}*X}}$. It is shown in
\cite{BH/NullityLoopcGraphsDM/10} that $d_{\mathcal{M}_{A}*X} =
\nu(A[X])$ for any symmetric or skew-symmetric matrix $A$.
Hence $d_{\mathcal{M}_{G}*X} = \nu(G[X])$.

If $u$ is a looped vertex, then $p_1(G)(y) = p_1(G\vertexrem u)(y) -
p_1(G*u\vertexrem u)(y)$ follows from Corollary~\ref{cor:rec_p1}
and the fact that $G*u$ is defined when $u$ is looped.

If $\{u,v\}$ is an edge where both $u$ and $v$ are not looped,
then by Corollary~\ref{cor:rec_p1} $p_1(G)(y) = p_1(\mathcal{M}_G)(y)
= p_1(\mathcal{M}_G\vertexrem u)(y) -
p_1(\mathcal{M}_G*u\vertexrem u)(y) = p_1(\mathcal{M}_G\vertexrem
u)(y) + p_1(\mathcal{M}_G*u*v\vertexrem u)(y)$. Now
$\mathcal{M}_G*u*v = \mathcal{M}_G*\{u,v\} =
\mathcal{M}_{G*\{u,v\}}$ where the second equality holds as
$G*\{u,v\}$ is defined. Moreover, $\mathcal{M}_G\vertexrem u =
\mathcal{M}_{G\vertexrem u}$ and
$\mathcal{M}_{G*\{u,v\}}\vertexrem u =
\mathcal{M}_{G*\{u,v\}\vertexrem u}$. Consequently,
$p_1(\mathcal{M}_G\vertexrem u)(y) +
p_1(\mathcal{M}_G*u*v\vertexrem u)(y) = p_1(G\vertexrem u)(y) +
p_1(G*\{u,v\}\vertexrem u)(y)$.

If $u$ is an isolated vertex of $G$, then $p_1(G)(y) =
(1-y) p_1(G\vertexrem u)(y)$ follows from
Corollary~\ref{cor:rec_p1}.
\end{Proof}

An example of the recursive computation is given in
Figure~\ref{fig:ex-recursive-pol} (ignore the caption of the
figure for now). The graph operations ${}*u$ on a looped vertex
and ${}*\{u,v\}$ on an unlooped edge are known as \emph{local
complementation} and \emph{edge local complementation},
respectively. Local complementation ``complements'' the edges
in the neighbourhood $N_G(u) = \{ v \in V \mid \{u,v\} \in
E(G), u \not= v \}$ of $u$ in $G$: for each $v,w \in N_G(u)$,
$\{v,w\}\in E(G)$ iff $\{v,w\} \not\in E(G*\{u\})$, and
$\{v\}\in E(G)$ iff $\{v\} \not\in E(G*\{u\})$ (the case
$v=w$). The other edges are left unchanged. We will not recall
the explicit graph theoretical definition of edge local
complementation in this paper. It can be found in, e.g.,
\cite{BH/PivotLoopCompl/09}.

\subsection{The Penrose Polynomial as a specialization of $Q(M)$}

Results for $p_1(M)(y) = P_{M+V*V}(y)$ can be straightforwardly translated to $P_{M}(y)$. Similar as for $p_1$, we have $P_M(1) = 0$ if $|V| > 0$. Since $p_1(M*Z)(y) = (-1)^{|Z|} p_1(M)(y)$, we have $P_{M+Z} = (-1)^{|Z|} P_{M}$ (this is also be easily verified from the definition of $P_{M}$).

\begin{Example}\label{ex:fano-penrose}
For the Fano matroid $F_7$ and and its dual $F_7^*$ we compute
the Penrose polynomial, cf.\ \cite{AignerPenroseBinMat}. Thus,
$P_{F_7}(y) = p_1(F_7*V+V)(y) = -p_1(F_7\dual V)(y) = -p_1(F_7)(y) \allowbreak = y^4
-8 y^3 +35 y^2 - 56 y + 28$, see Example~\ref{ex:fano-p1}. Note
that $F_7 \dual V = F_7$, so $F_7^*  +V = F_7^*$ and
$P_{F_7^*} = 0$.
\end{Example}

If $M$ is an equicardinal set system, then $\max(M*V) = M*V$ and thus, for every $X \subseteq V$, the sets of $M*V$ are also present in $M*V\dual X$. Consequently, $d_{M*V\dual X} \le d_{M*V}$. Therefore, the degree of $P_M(y)$ can be at most $d_{M*V}$. In fact, it cannot be less than that value, without becoming nontrivial, generalizing \cite[Proposition~2]{AignerPenroseBinMat}.
\begin{theorem}\label{thm:zero-penrose}
Let $M$ be an equicardinal set system over $V$. Then $P_M(y) = 0$ iff the degree of $P_M(y)$ is smaller than $d_{M*V}$ iff $M+X=M$ for some $X \subseteq V$ with $|X|$ odd.
\end{theorem}
\begin{Proof}
Trivially, if $P_M(y) = 0$, then the degree $d$ of $P_M(y)$ is
smaller than $d_{M*V}$. Assume $d<d_{M*V}$. Since $X =
\emptyset$ in $d_{M*V\dual X}$ contributes $+1$ to the
coefficient of $y^{d_{M*V\dual X}}$, there must be an $X \subseteq V$ with
$|X|$ odd contributing $-1$ to the coefficient of $y^{d_{M*V\dual X}}$. Since
all sets of $M*V$ have equal cardinality, $d_{M*V\dual X} =
d_{M*V}$ implies $M*V\dual X = M*V$.
Hence, $M+X=M$.
Finally, if $M+X=M$ for some $X$ of odd cardinality, then $P_M(y)=0$ since $P_{M}(y) = (-1)^{|X|} P_{M+X}(y)$.
\end{Proof}
Note that, for $v \in V$, $M+v = M$ iff $v$ is a coloop of $M$. In \cite{AignerPenroseBinMat} a necessary condition is given of the binary matroids $M$ for which $P_M(y) = 0$: either $M$ has a coloop, or $M$ has $F^*_7$ as a minor. Note that the uniform matroid $U_{3,5}$ (which is vf-safe but not binary) we also have $P_M(U_{3,5})(y) = 0$, as $U_{3,5} + V = U_{3,5}$ is readily verified.

We now carry over properties from $p_1(M)(y)$ to the Penrose polynomial $P_M(y)$. We find recursive relations that characterize the Penrose polynomial. It seems such recursive relations do not exist for $4$-regular graphs or matroids---one needs to step out of those domains, like in \cite{EllisMonaghan/PenroseEmb/2011}, where graphs embedded in surfaces are considered.

We first show that $P_M(y)$ fulfills the following characterizing recursive
relation.
\begin{theorem}\label{thm:recursive-p2}
Let $M$ be a vf-safe \dmatroid, and let $u\in V$. If $u$ is nonsingular in $M\dual V$, then $P_M(y) = P_{M*u\vertexrem u}(y) -  P_{M \dual u \vertexrem u}(y)$. If $u$ is a coloop of $M\dual V$, then $P_M(y) = (1-y) P_{M*u\vertexrem u}(y)$, and if $u$ is a loop of $M\dual V$, $P_M(y) = (y-1) P_{M \dual u \vertexrem u}(y)$. Finally, if $V = \emptyset$, then $P_M(y) = 1$.
\end{theorem}
\begin{Proof}
Let $u$ is nonsingular in $M\dual V$ and let $V'= V\vertexrem
\{u\}$. We have $P_M(y) = p_1(M*V+V)(y) = p_1(M*V+V\vertexrem u)(y) -
p_1(M*V+V*u\vertexrem u)(y)  = p_1(M*u+u \vertexrem u*V'+V')(y) -
p_1(M*u+u*u\vertexrem u *V'+V')(y) = P_{M*u+u\vertexrem u}(y) -  P_{M \dual u \vertexrem u}(y)$. Finally observe that $N+u\vertexrem u =
N\vertexrem u$ for any set system $N$. The cases where $u$ is singular in $M\dual
V$ is proved similarly.
\end{Proof}
Note that Theorem~\ref{thm:recursive-p2} holds for every set system $M$ such that $M\dual V$ is a \dmatroid. For notational convenience we assume that $M$ is vf-safe, but the reader may easily recover the loss of generality.

We now show that, for a binary matroid $M$, $P_M(y)$ may also be viewed (up to a sign) as the graph polynomial $p_1(G)(y)$ where $G$ is a fundamental graph of $M$.
\begin{theorem} \label{thm:penrose_graphpol}
Let $M$ be a binary matroid and $Z$ a basis of $M$. Then $M
\dual V*Z = \mathcal{M}_{G}$ for some graph $G$ and $P_M(y) =
(-1)^{\nu(M)} p_1(G)(y)$.
\end{theorem}
\begin{Proof}
Since $Z \in M = \max(M)$, we have $Z \in M \dual V$, and so
$\emptyset \in M \dual V*Z$. Hence $M \dual V*Z =
\mathcal{M}_{G}$ for some graph $G$. Moreover, $P_M(y) =
P_M(y) = P_{\mathcal{M}_{G}*Z\dual V}(y) =
p_1(\mathcal{M}_{G}*Z\dual V*V+V)(y) =
p_1(\mathcal{M}_{G}*(V\setminus Z))(y) = (-1)^{|V\setminus
Z|}p_1(\mathcal{M}_{G})(y)$. As $|V\setminus Z| = \nu(M)$, the
result follows.
\end{Proof}

\begin{Example} \label{ex:recursive-diamond}
Consider the cycle matroid $M_D$ of the diamond graph $D$ with
edge set $V=\{1,2,\dots,5\}$, cf.\ Figure~\ref{fig:diamond}. Its
eight bases, the spanning trees of $D$, are all subsets of $V$
of cardinality $3$ except for the forbidden triangles
$\{1,4,5\}$ and $\{2,3,5\}$.
We compute the Penrose polynomial $P_{M_D}(y)$ recursively with
the help of Theorem~\ref{thm:penrose_graphpol} and
Lemma~\ref{lem:p1_graph}.
\begin{figure}\sidecaption
\unitlength 0.8mm%
\centerline{%
\begin{picture}(24,30)(-2,-5)
\gasset{AHnb=0,Nw=2.5,Nh=2.5,Nframe=n,Nfill=y}
  \node(1)(00,20){}
  \node(2)(00,00){}
  \node(3)(20,00){}
  \node(4)(20,20){}
  \drawedge(4,2){5}
  \drawedge(2,1){4}
  \drawedge(3,2){3}
  \drawedge(4,3){2}
  \drawedge(1,4){1}
\end{picture}%
}%
\caption{Diamond graph $D$.}\label{fig:diamond}
\end{figure}
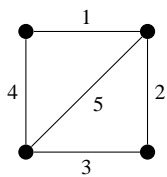

We have that $M_D\dual V$ is obtained from $M_D$ by adding the
family of sets $\{
\{1\},\{2\},\allowbreak\{3\},\{4\},\{1,2\},\{1,3\},\{2,4\},\{3,4\} \}$.
Now, binary \dmatroid $M_D\dual V * \{1,2,3\}$ contains
$\emptyset$, which hence represents a graph $G$.
%
The sets in $M_D\dual V * \{1,2,3\}$ of cardinality one or two
(which uniquely determine $G$) are $\{
\{2\},\{3\},$
$\{1,2\},\{1,3\},\{1,4\},\{2,3\},\allowbreak
\{2,4\},\{2,5\},\allowbreak \{3,4\},\{3,5\} \}$.
%

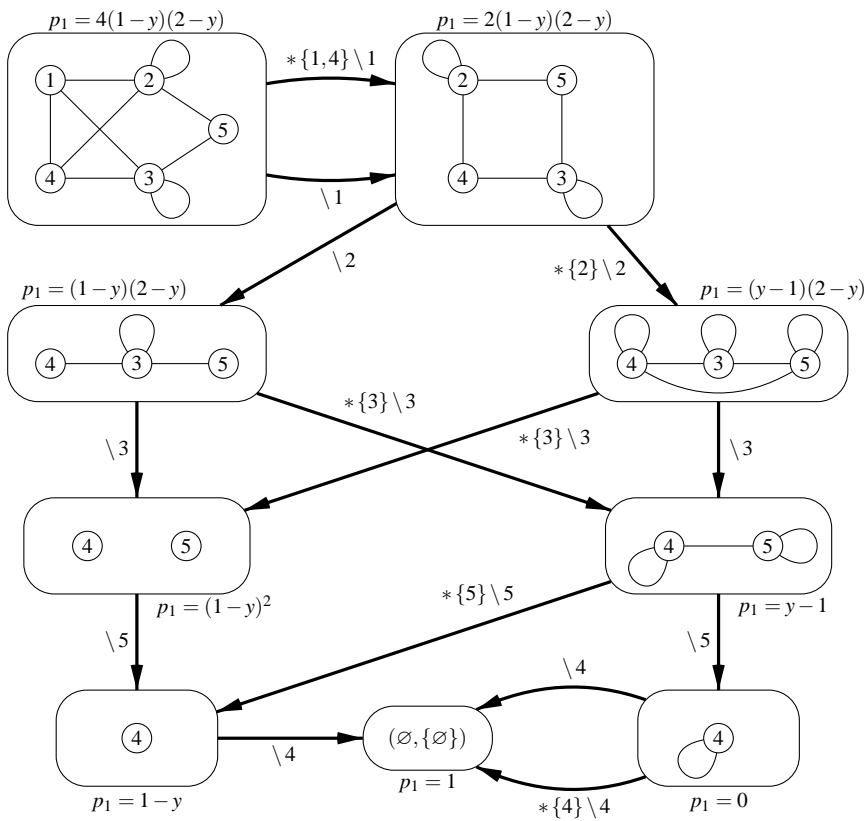
\begin{figure}[!t]
\begin{center}
\unitlength 0.85mm
\newcommand{\UnitPict}{\unitlength0.65mm}
{%
\begin{picture}(130,130)(0,70)
\node[Nw=40,Nh=30,Nframe=y](G)(20,180){%
{ \UnitPict
  \begin{picture}(45,30)
  \gasset{AHnb=0,Nw=6,Nh=6,Nframe=y,Nfill=n}
    \node(1)(05,25){$1$}
    \node(2)(25,25){$2$}
    \node(3)(25,05){$3$}
    \node(4)(05,05){$4$}
    \node(5)(40,15){$5$}
  \drawedge(1,4){}
  \drawedge(1,2){}
  \drawedge(1,3){}
  \drawedge(2,4){}
  \drawedge(3,4){}
  \drawedge(2,5){}
  \drawedge(3,5){}
  \drawloop[loopangle=45,loopdiam=7](2){}
  \drawloop[loopangle=-45,loopdiam=7](3){}
\end{picture}}
}
  \node[Nadjust=wh,Nframe=n](P)(20,197){$p_1 = 4(1-y)(2-y)$}
\node[Nw=40,Nh=30,Nframe=y](G1)(80,180){%
{ \UnitPict
  \begin{picture}(35,30)
  \gasset{AHnb=0,Nw=6,Nh=6,Nframe=y,Nfill=n}
    \node(2)(05,25){$2$}
    \node(3)(25,05){$3$}
    \node(4)(05,05){$4$}
    \node(5)(25,25){$5$}
  \drawedge(2,4){}
  \drawedge(4,3){}
  \drawedge(3,5){}
  \drawedge(5,2){}
  \drawloop[loopangle=135,loopdiam=7](2){}
  \drawloop[loopangle=-45,loopdiam=7](3){}
\end{picture}}
}
  \node[Nadjust=wh,Nframe=n](P)(80,197){$p_1 = 2(1-y)(2-y)$}
\node[Nw=40,Nh=15,Nframe=y](G12)(20,145){%
{ \UnitPict
  \begin{picture}(35,10)
  \gasset{AHnb=0,Nw=6,Nh=6,Nframe=y,Nfill=n}
    \node(3)(17.5,03){$3$}
    \node(4)(0,03){$4$}
    \node(5)(35,03){$5$}
  \drawedge(5,3){}
  \drawedge(3,4){}
  \drawloop[loopangle=90,loopdiam=7](3){}
\end{picture}}
}
  \node[Nadjust=wh,Nframe=n](P)(15,155){$p_1 = (1-y)(2-y)$}
\node[Nw=40,Nh=15,Nframe=y](G12*)(110,145){%
{ \UnitPict
  \begin{picture}(35,10)
  \gasset{AHnb=0,Nw=6,Nh=6,Nframe=y,Nfill=n}
    \node(3)(17.5,03){$3$}
    \node(4)(00,03){$4$}
    \node(5)(35,03){$5$}
  \drawedge[curvedepth=-6](4,5){}
  \drawedge(5,3){}
  \drawedge(3,4){}
  \drawloop[loopangle=90,loopdiam=7](3){}
  \drawloop[loopangle=90,loopdiam=7](4){}
  \drawloop[loopangle=90,loopdiam=7](5){}
\end{picture}}
}
  \node[Nadjust=wh,Nframe=n](P)(120,155){$p_1 = (y-1)(2-y)$}
\node[Nw=35,Nh=15,Nframe=y](G-45L)(20,115){%
{ \UnitPict
  \begin{picture}(30,10)
  \gasset{AHnb=0,Nw=6,Nh=6,Nframe=y,Nfill=n}
    \node(4)(05,05){$4$}
    \node(5)(25,05){$5$}
\end{picture}}
}
  \node[Nadjust=wh,Nframe=n](P)(32,105.50){$p_1 = (1-y)^2$}
\node[Nw=35,Nh=15,Nframe=y](G-45R)(110,115){%
{ \UnitPict
  \begin{picture}(30,10)
  \gasset{AHnb=0,Nw=6,Nh=6,Nframe=y,Nfill=n}
    \node(4)(05,05){$4$}
    \node(5)(25,05){$5$}
  \drawedge(4,5){}
  \drawloop[loopangle=225,loopdiam=7](4){}
  \drawloop[loopangle=0,loopdiam=7](5){}
\end{picture}}
}
  \node[Nadjust=wh,Nframe=n](P)(120,105.5){$p_1 = y-1$}
\node[Nw=25,Nh=15,Nframe=y](G-4L)(20,85){%
{ \UnitPict
  \begin{picture}(30,10)
  \gasset{AHnb=0,Nw=6,Nh=6,Nframe=y,Nfill=n}
    \node(4)(15,05){$4$}
\end{picture}}
}
  \node[Nadjust=wh,Nframe=n](P)(20,75){$p_1 = 1-y$}
\node[Nw=25,Nh=15,Nframe=y](G-4R)(110,85){%
{ \UnitPict
  \begin{picture}(30,10)
  \gasset{AHnb=0,Nw=6,Nh=6,Nframe=y,Nfill=n}
    \node(4)(15,05){$4$}
  \drawloop[loopangle=225,loopdiam=7](4){}
\end{picture}}
}
  \node[Nadjust=wh,Nframe=n](P)(110,75){$p_1 = 0$}
\node[Nw=20,Nh=10,Nframe=y](G0)(65,85){$(\varnothing,\{\varnothing\})$}
  \node[Nadjust=wh,Nframe=n](P)(65,78){$p_1 = 1$}
    \drawedge[curvedepth=-8,ELside=r,linewidth=0.5,AHLength=4.0,AHlength=4,AHangle=14](G,G1){${}\vertexrem 1$}
    \drawedge[curvedepth=8,linewidth=0.5,AHLength=4.0,AHlength=4,AHangle=14](G,G1){${}*\{1,4\}\vertexrem 1$}
    \drawedge[linewidth=0.5,AHLength=4.0,AHlength=4,AHangle=14](G1,G12){${}\vertexrem 2$}
    \drawedge[ELside=r,linewidth=0.5,AHLength=4.0,AHlength=4,AHangle=14](G1,G12*){${}*\{2\}\vertexrem 2$}
    \drawedge[ELside=r,linewidth=0.5,AHLength=4.0,AHlength=4,AHangle=14](G12,G-45L){${}\vertexrem 3$}
    \drawedge[ELpos=40,linewidth=0.5,AHLength=4.0,AHlength=4,AHangle=14](G12,G-45R){${}*\{3\}\vertexrem 3$}
    \drawedge[ELpos=30,linewidth=0.5,AHLength=4.0,AHlength=4,AHangle=14](G12*,G-45L){${}*\{3\}\vertexrem 3$}
    \drawedge[linewidth=0.5,AHLength=4.0,AHlength=4,AHangle=14](G12*,G-45R){${}\vertexrem 3$}
    \drawedge[ELside=r,linewidth=0.5,AHLength=4.0,AHlength=4,AHangle=14](G-45L,G-4L){${}\vertexrem 5$}
    \drawedge[ELside=r,ELpos=40,linewidth=0.5,AHLength=4.0,AHlength=4,AHangle=14](G-45R,G-4L){${}*\{5\}\vertexrem 5$}
    \drawedge[ELside=r,linewidth=0.5,AHLength=4.0,AHlength=4,AHangle=14](G-45R,G-4R){${}\vertexrem 5$}
    \drawedge[ELside=r,linewidth=0.5,AHLength=4.0,AHlength=4,AHangle=14](G-4L,G0){${}\vertexrem 4$}
    \drawedge[curvedepth=8,linewidth=0.5,AHLength=4.0,AHlength=4,AHangle=14](G-4R,G0){${}*\{4\}\vertexrem 4$}
    \drawedge[ELside=r,curvedepth=-8,linewidth=0.5,AHLength=4.0,AHlength=4,AHangle=14](G-4R,G0){${}\vertexrem 4$}
\end{picture}
}
\end{center}
\caption{Recursive computation of the Penrose polynomial $P_{M_D}(y)$, cf.\ Example~\ref{ex:recursive-diamond}.
With each graph $F$ we give the polynomial $p_1 = p_1(F)(y)$.}\label{fig:ex-recursive-pol}
\end{figure}
By applying the graph operations of local complementation and
edge local complementation, we determine the polynomial
$p_1(G)(y) = 4(1-y)(2-y)$, see
Figure~\ref{fig:ex-recursive-pol}. Thus $P_{M_D}(y) =
4(y-1)(2-y)$.
\end{Example}

We now extend the result from
\cite[Proposition~1]{AignerPenroseBinMat} that for an Eulerian
matroid $M$ the value $P_M(2)$ of the Penrose polynomial of a
binary matroid equals the size of its cocycle space. We also
consider $P_M(-1)$.

\begin{theorem} \label{thm:p2_at_2}
Let $M$ be a set system such that $M\dual V$ is even. Then
$P_M(y) = (-1)^{d_{M*V}}\allowbreak Q_{[1,1,0]}(M\dual V)(-y)$. In particular, $P_M(-1) = (-1)^{d_{M*V}}2^{|V|}$. If $M$ is moreover a vf-safe \dmatroid, then $P_M(2) = \allowbreak
(-1)^{d_{M}+d_{M*V}+|V|} \allowbreak 2^{d_{M}}$.
\end{theorem}
\begin{Proof}
Note that $M\dual V$ is an even set system iff $M*V+V$ has that property. Thus, we have by Lemma~\ref{lem:even-p1-q1} $P_M(y) = p_1(M*V+V)(y) = (-1)^{d_{M*V+V}}Q_{[1,1,0]}(M*V+V)(-y) = (-1)^{d_{M*V}}Q_{[1,1,0]}(M\dual V)(-y)$, where in the last equality we use that $d_{N+V} = d_N$ and $Q_{[1,1,0]}(N*V)(y) = Q_{[1,1,0]}(N)(y)$ for every set system $N$ over $V$. Note also that $Q_{[1,1,0]}(N)(1) = 2^{|V|}$ for every set system $N$ over $V$.

By \cite{BH/InterlacePolyDM/14} we have $Q_{[1,1,0]}(N)(-2) =
(-1)^{|V|} (-2)^{d_{N\dual V}}$ for any vf-safe \dmatroid
$N$. Thus, $P_M(2) = (-1)^{d_{M*V}} \allowbreak Q_{[1,1,0]}(M\dual V)(-2) \allowbreak =
(-1)^{d_{M*V}+|V|} (-2)^{d_{M}}$.
\end{Proof}



As another example, for the nonbinary matroid $U_{2,5}$ we have $U_{2,5} \dual V = U_{2,5}$, so $U_{2,5} \dual V$ is even. Then $P_{U_{2,5}}(y) = p_1(U_{2,5} *V+V)(y) = p_1(U_{2,5}\dual V *V)(y) = p_1(U_{2,5}*V)(y) = p_1(U_{3,5} )(y) = y^3 - 5y^2 + 10 y - 10 + 5 y -y^2 = y^3 -6y^2 +15y -10$. Thus $P_{U_{2,5}}(2) = 4 = 2^{\rho(U_{2,5})}$, as predicted by Theorem~\ref{thm:p2_at_2}.

We now turn to binary matroids. By
Theorem~\ref{thm:p2_at_2} and Theorem~\ref{thm:bipartite_char}(2) we have the following
result (we also use that $Q_{[1,1,0]}(N)(1) = 2^{|V|}$ for every set system $N$).
\begin{corollary} \label{cor:eval_penrose}
Let $M$ be an Eulerian binary matroid. Then we have $P_M(y) =
(-1)^{\nu(M)} \allowbreak Q_{[1,1,0]}(M\dual V)(-y)$, $P_M(-1) =
(-1)^{\nu(M)}2^{|V|}$, and $P_M(2) = 2^{\rho(M)}$.
\end{corollary}



We now generalize the equality between the Penrose polynomial
at $-2$ and the Tutte polynomial at $(0,-3)$, given in
\cite[Theorem~2]{AignerPenroseBinMat}, from binary matroids to
vf-safe matroids. It is well known that $|t_M(0,-3)|$ is equal to the number of nowhere-zero $4$-flows of a binary matroid $M$.
\begin{corollary} \label{cor:tutte-at-0-3}
Let $M$ be a vf-safe matroid. Then $P_M(-2) = 2^{\rho(M)}
t_M(0,-3)$ where $t_M(x,y)$ is the Tutte polynomial.
\end{corollary}
\begin{Proof}
First we need the following auxiliary result. Let $a,b,c,d$ be arbitrary values. Then $Q_{[a,b,c]}(M)(-2) = Q_{[a+d,b+d,c+d]}(M)(-2)$ for vf-safe \dmatroids $M$. This equality is a special case of \cite[Theorem~7]{BH/InterlacePolyDM/14}, which is more generally stated there in terms of tight multimatroids (the conversion to vf-safe \dmatroids is similar as in the proof of Theorem~38.2 in \cite{BH/InterlacePolyDM/14}).

By the above auxiliary result, we have $P_M(-2)\allowbreak = Q_{[0,1,-1]}(M)(-2) = Q_{[1,2,0]}(M)(-2)$. Now by Proposition~\ref{prop:transtition-tutte}, we obtain $Q_{[1,2,0]}\allowbreak(M)\allowbreak(-2)\allowbreak = 2^{\rho(M)} t_M(0, -3)$.
\end{Proof}

\bibliography{../geneassembly}

\end{document}